\documentclass[11pt]{article}

\usepackage{amsmath,amsthm,amsfonts}

\usepackage{amssymb}

\usepackage{epsfig}

\usepackage{rotating}

\usepackage{subfigure}

\newcommand{\p}[2]{\frac{\partial#1}{\partial#2}}

\def\hR{{\hat{R}}}
\def\eps{\epsilon}

\begin{document}

\title{Dimensional reduction as a tool for mesh refinement and tracking singularities of PDEs}
\author{Panagiotis Stinis \\ 
Department of Mathematics \\
Lawrence Berkeley National Laboratory \\
    Berkeley, CA 94720 \\
   stinis@math.lbl.gov} 
\date {}

\maketitle

\begin{abstract}
We present a collection of algorithms which utilize dimensional reduction to perform mesh refinement and study possibly singular solutions of time-dependent partial differential equations. The algorithms are inspired by constructions used in statistical mechanics to evaluate the properties of a system near a critical point. The first algorithm allows the accurate determination of the time of occurrence of a possible singularity. The second algorithm is an adaptive mesh refinement scheme which can be used to approach efficiently the possible singularity. Finally, the third algorithm uses the second algorithm until the available resolution is exhausted (as we approach the possible singularity) and then switches to a dimensionally reduced model which, when accurate, can follow faithfully the solution beyond the time of occurrence of the purported singularity. An accurate dimensionally reduced model should dissipate energy at the right rate.

We construct two variants of each algorithm. The first variant assumes that we have actual knowledge of the reduced model. The second variant assumes that we know the form of the reduced model, i.e. the terms appearing in the reduced model, but not necessarily their coefficients. In this case, we also provide a way of determining the coefficients. We present numerical results for the Burgers equation with zero and nonzero viscosity to illustrate the use of the algorithms.  
\end{abstract}


\section*{Introduction}
The problem of constructing dimensionally reduced models for large systems of ordinary differential equations (this covers the case of partial differential equations after discretization or series expansion of the solution) has received considerable attention in the last decade e.g. see the review papers \cite{givon, CS05}. The construction of an accurate reduced model has advantages beyond the obvious one of predicting the correct behavior for a reduced set of variables. For example, one can use the reduced model to perform more efficiently tasks like e.g. filtering \cite{GSW}, which can be prohibitively expensive for the original (full-dimensional) system. Following this line of thought, we present here a collection of algorithms that are based on dimensional reduction and which can be used to perform mesh refinement and investigate possibly singular solutions of partial differential equations. The algorithms are inspired by constructions used in statistical mechanics to evaluate the properties of a system near a critical point \cite{goldenfeld}.

The algorithms we present address three objectives: i) they provide a way of accurately monitoring the progress of a simulation towards underresolution, thus providing us as a byproduct with the time of occurrence of the possible singularity; ii) they allow the formulation of a mesh refinement scheme and iii) they allow the simulation to follow the purported singularity beyond the point where the resolution is exhausted by switching to a dimensionally reduced model which dissipates energy. The end result of the algorithms is that we are able to reach the time of interesting dynamics of the equation much more efficiently compared to an algorithm that simply starts with the maximum available resolution.

The task of investigating numerically the appearance of a singularity is subtle. Clearly, since all calculations are performed with finite resolution and a singularity involves an infinity of active scales we can only come as close to the singularity as our resolution will allow. After that point, either we stop our calculations and conclude that a singularity may be present close to the time instant we stopped or we switch, if available, to a model that drains energy at the correct rate out of the set of resolved variables. We emphasize here that, up to some time, the evolution towards a near-singular solution can be identical to the evolution towards a singular solution. If we do not have enough resolution to go beyond the time instant after which the two evolutions start deviating, we cannot claim with certainty the presence of a singularity. In other words, given adequate resolution we can eliminate the possibility of a singularity but it may be very hard, if not impossible, to prove through a finite calculation, no matter how large, that a singularity exists. We will come back to this point in the discussion of the algorithms in Section \ref{algo} and through numerical examples in Sections \ref{numerics} and \ref{numerics2}.

The paper is organized as follows. In Section \ref{algo} we present the ideas behind the construction of the three algorithms.  There are two variants of the algorithms. The algorithms for the case when we have knowledge of the functional form of the reduced model and its related coefficients are presented in Section \ref{know}. The algorithms for the case when we only know the functional form of the reduced model but not necessarily its coefficients are presented in Section \ref{know2}. Each of the Sections \ref{know} and \ref{know2} contains numerical results for the Burgers equation with and without viscosity. Section \ref{conclusions} contains some conclusions.


\section{The main construction}\label{algo}

Suppose that we are interested in the possible development of singularities in the solution $v(x,t)$ of a partial differential equation (PDE)
$$ v_t + H (t,x,v,v_x,...)=0 $$
where $H$ is a, in general nonlinear, operator and $x \in D \subseteq \mathbb{R}^d$ (the constructions extend readily to the case of systems of partial differential equations). After spatial discretization or expansion of the solution in series, the PDE transforms into a system of ordinary differential equations (ODEs). For simplicity we restrict ourselves to the case of periodic boundary conditions, so that a Fourier expansion of the solution leads to system of ODEs for the Fourier coefficients. To simulate the system for the Fourier coefficients we need to truncate at some point the Fourier 
expansion. Let $F \cup G$ denote the set of Fourier modes retained in the series, where we have split the Fourier modes in two sets, $F$ and $G.$ 
We call the modes in $F$ resolved and the modes in $G$ unresolved. One can construct, in principle, an 
exact reduced model for the modes in $F$ e.g. through the Mori-Zwanzig formalism \cite{CHK00} (we do not deal here with the complications of constructing a good reduced model).

The main idea behind the algorithms is that the evolution of moments of the reduced set of modes, for example $L_p$ norms of the modes in $F$, should be the same whether computed from the full or the reduced system. This is a generalization to time-dependent systems of the principle used in the theory of equilibrium phase transitions to compute the critical exponents \cite{goldenfeld,S05}. The idea underlying the computation of the critical exponents is that while the form of the reduced system equations is important, one can extract even more information by looking at how the form of the reduced system is related to the form of the original (full dimensional) system. This is expressed mathematically through the chain rule and leads to the formation of a matrix (known as the renormalization matrix in the statistical mechanics literature) whose eigenvalues allow the determination of the critical exponents. In the case of an equilibrium phase transition, one can use this approach to compute properties of the system near the critical value of the governing parameter (e.g. near the critical temperature). If one thinks of time as a parameter and of the time of occurrence of a singularity as a critical value of the parameter, then by looking at the eigenvalues of an appropriately defined renormalization matrix, we can study the behavior of a system near the time of formation of a singularity. We caution the reader that even though our motivation for the present construction came from the theory of equilibrium phase transitions, we do not advocate that a singularity is a phase transition in the conventional sense. It can be thought of as a transition from a strong solution to an appropriately defined weak solution but one does not have to push the analogy further. We want to point here that the problem we are addressing is different from the subject known as dynamic critical phenomena (see Ch. 8 in \cite{goldenfeld}). There, one is interested in the computation of time-dependent quantities as a controlling parameter, other than time, reaches its critical value. In our case, time {\it is} the controlling parameter and we are interested in the behavior of the solution as time reaches a critical value. 

The above arguments can be made more precise. The original system of equations for the modes $F \cup G$ is given by 
$$\frac{du(t)}{dt} = R (t,u(t)),$$
where $u = ( \{u_k\}), \; k \in F \cup G$ is the vector of Fourier coefficients of $u$ and $R$ is the Fourier transform of the operator $H.$ The system should be supplemented with an initial condition $u(0)=u_0.$ The vector of Fourier coefficients can be written as $ u = (\hat{u}, 
\tilde{u}),$ where $ \hat{u}$ are the resolved modes (those in $F$) and $\tilde{u}$ the unresolved ones (those in $G$). Similarly, for the right hand sides (RHS) we have $R(t,u) = (\hat{R}(t,u), \tilde{R}(t,u)).$ Note that the RHS of the resolved modes involves both resolved and unresolved modes. In anticipation of the construction of a reduced model we can rewrite the RHS as $R(t,u)=R^{(0)}(t,u) = (\hat{R}^{(0)}(t,u), \tilde{R}^{(0)}(t,u)).$ Furthermore, we can decompose $\hat{R}^{(0)}(t,u)$ as 
$$\hat{R}^{(0)}(t,u(t)) = \sum_{i=1}^{m} a^{(0)}_i \hat{R}^{(0)}_i(t,u(t)).$$
The equations for the reduced set of modes in $F$ is written as 
\begin{equation}\label{full}
\frac{d\hat{u}(t)}{dt} = \hat{R}^{(0)} (t,u(t))=\sum_{i=1}^{m} a^{(0)}_i \hat{R}^{(0)}_i (t,u(t))
\end{equation}
Note that not all the coefficients $a^{(0)}_i, \; i =1,\ldots,m$ have to be nonzero. As is standard in renormalization theory \cite{binney}, one augments (with zero coefficients) the RHS of the equations in the full system by terms whose form is the same as the terms appearing in the RHS of the equations for the reduced model. Dimensional 
reduction transforms the vector $a^{(0)}=(a^{(0)}_1,\ldots,a^{(0)}_m)$ to $a^{(1)}=(a^{(1)}_1,\ldots,a^{(1)}_m).$ The reduced model for the modes in $F$ is given by 
\begin{equation}\label{reduced}
\frac{d\hat{u}'(t)}{dt} = \hat{R}^{(1)} (t,\hat{u}'(t))=\sum_{i=1}^{m} a^{(1)}_i \hat{R}^{(1)}_i (t,\hat{u}'(t))
\end{equation}
with initial condition $\hat{u}'(0)=\hat{u}_0.$ We emphasize that the functions $\hat{R}^{(1)}_i, \; i=1,\ldots,m$ have the same form as the functions $\hat{R}^{(0)}_i, \; i=1,\ldots,m,$ but they are defined {\it only} on the reduced set of modes $F.$ This allows one to determine the relation of the full to the reduced system by focusing on the change of the vector $a^{(0)}$ to $a^{(1)}.$ Also, the vectors $a^{(0)}$ and $a^{(1)}$ do not have to be constant in time. This does not change the analysis that follows.

Define $m$ quantities 
$\hat{E}_i, \; i=1,\ldots,m$ involving only modes in $F.$ For example, these could be $L_p$ norms of the reduced set of modes. To proceed we require that these quantities' rates of change are the same when computed from (\ref{full}) and (\ref{reduced}), i.e. 
\begin{equation}\label{conditions}
\frac{d\hat{E}_i(\hat{u})}{dt} = \frac{d\hat{E}_i(\hat{u}')}{dt}, \; i=1,\ldots,m.
\end{equation}
Note that similar conditions, albeit time-independent, lie at the heart of the renormalization group theory for equilibrium systems (\cite{binney} p. 154). In fact, it is these conditions that allow the definition and calculation of the (renormalization) matrix whose eigenvalues are used to calculate the critical exponents. In the current (time-dependent) setting, the renormalization matrix is defined by 
differentiating $\frac{d\hat{E}_i(\hat{u})}{dt}$ with respect to $a^{(0)}$ and using (\ref{conditions}) to 
obtain
\begin{equation}\label{rngmatrix1}
\p{}{a^{(0)}_j}\biggl(\frac{d\hat{E}_i(\hat{u})}{dt}\biggr)=\sum_{k=1}^{m}\p{}{a^{(1)}_k}\biggl(\frac{d\hat{E}_i(\hat{u}')}{dt}\biggr) \p{a^{(1)}_k}{a^{(0)}_j}, \; i,j=1,\ldots,m.
\end{equation}
We define the renormalization matrix $M_{kj}= \p{a^{(1)}_k}{a^{(0)}_j}, \; k,j=1,\ldots,m,$ as well as the matrices $A_{kj}=\p{}{a^{(0)}_j}\biggl(\frac{d\hat{E}_k(\hat{u})}{dt}\biggr), \; k,j=1,\ldots,m$ and $B_{kj}=\p{}{a^{(1)}_j}\biggl(\frac{d\hat{E}_k(\hat{u}')}{dt}\biggr), \; k,j=1,\ldots,m .$ Equations (\ref{rngmatrix1}) can 
be written in matrix form as
\begin{equation}\label{rngmatrix2}
A=MB
\end{equation} 
The eigenvalues of the matrix $M$ contain information about the behavior of the reduced system relative to the full system. In fact, the different eigenvalues provide information as to how the different terms appearing in the RHS of the equations for the full and reduced system contribute to the evolution of the quantities $\hat{E}_i, \; i=1,\ldots,m.$ Thus, if we monitor the behavior of the eigenvalues as a function of time, we should 
be able to determine whether the two systems are coming closer or deviating. We know that while 
the full system can become underresolved, an accurate reduced model can allow us to follow the 
evolution of the reduced set of modes $F$ after the full system has become underresolved. As a consequence, the behavior of the reduced system should deviate from that of the full system as 
the full system becomes underresolved. Our algorithms exploit this observation to perform the following tasks: i) detect when the full system becomes underresolved, ii) increase the resolution as needed and iii) switch to a reduced model when no more refinement is possible.

Before we present the algorithms, there is a seemingly paradoxical situation that needs to be 
resolved. We have made the statement that a reduced model coming from the full system can 
continue to describe accurately the behavior of the reduced set of modes after the full system has 
become underresolved. This means that a reduced system coming from an eventually underresolved full system remains well resolved. Furthermore, it means that the reduced system is a good model for the 
solution of the PDE for all times, while the full system is a good model for the solution of the PDE only for a finite time interval. Thus, after some time, the reduced system must become a bad model for the 
full system while remaining a good model for the solution of the PDE. This can happen if the reduced (renormalized) model captures the essential part of the dynamics of the PDE while disregarding the parts that eventually lead to the underresolution of the full (bare) system. This situation is not new in the physics literature. Work on equilibrium phase transitions and quantum field theory in the late 60s and early 70s, showed that bare theories that are only good down to a certain length scale, can give rise to renormalized theories (at larger length scales) that can provide excellent agreement with physical experiments (e.g. \cite{weinberg}, Ch. 14  in \cite{binney}, Ch. 9 in \cite{goldenfeld}). This is the result of the suppression, as we go to larger length scales, of the "bad" terms in the bare theory, i.e. the terms that limit the applicability of the bare theory. In the physics terminology the terms that survive the renormalization procedure are called relevant, while the ones that get weaker as we go to larger scales are called irrelevant. In our setting, the "physical world" is provided by the PDE, the bare theory by the full system and the renormalized model by the reduced model. Similarly, the restriction of computational power in numerical simulations is the analog of the restriction of the available energy regimes in particle physics experiments. During the time interval when the bare theory is still valid, i.e. when there is no activity in scales smaller than the ones described by the full system, a good renormalized (reduced) model which keeps the relevant terms and discards the irrelevant ones is a good model for the solution of the PDE at the scales where the renormalized model is supposed to hold. After this time interval, the full system breaks down as an approximation, but a reduced model that keeps only the relevant terms will not suffer from the same problem because it neglects the dangerous irrelevant terms. Thus, a reduced model can continue being a good model for the solution of the PDE after the full system has ceased to be valid.

\section{Case I: Knowledge of the terms and coefficients of the reduced model}\label{know}
We begin our presentation of the algorithms with the algorithms for the case when we have complete knowledge of the reduced model, i.e. knowledge of $\hR^{(1)}$ and $a^{(1)}.$

\subsection{How to locate the time of occurrence of a possible singularity}\label{location}
The task of computing the time of occurrence of a singularity with a finite calculation, i.e. a finite number of modes, is delicate. The reason is that a finite calculation will eventually run out of resolution as we approach the singularity. An underresolved computation is not trustworthy. Yet, we want to push our 
time of prediction as close to the singularity as possible. This situation is equivalent to predicting with a finite system size the critical value of a parameter that controls a phase transition. A phase transition is a phenomenon that requires an infinitely large system. Any finite calculation can suffer from finite-size effects which contaminate the prediction of the critical value of the controlling parameter. What one 
has to do is to monitor the degree of underresolution of the full system and decide when the underresolution dominates the results, namely when the useful informational content of our calculation shrinks to zero. We propose to do that by looking at the temporal behavior of the eigenvalues of the 
renormalization matrix $M$ defined in the previous section. To do that, one has to evolve simultaneously the systems (\ref{full}) and (\ref{reduced}), compute the matrices $A$ and $B,$ solve (\ref{rngmatrix2}) for $M$ and compute the eigenvalues of $M.$

The subtle point here is that the definition 
of the renormalization matrix rests on the assumption that the conditions (\ref{conditions}) hold. But, as the calculation becomes underresolved, these conditions will only hold approximately. Thus, one has 
to monitor the validity of the conditions (\ref{conditions}) in addition to monitoring the behavior of the 
eigenvalues of the matrix $M.$

We can make some observations about what the values and the temporal behavior of the eigenvalues should be. Inspection of the entries of the matrices $A$ and $B$ reveals that the entries of these matrices are exactly the contributions of the different terms appearing in the RHS of (\ref{full}) and (\ref{reduced}) respectively, to the rates of change of the quantities $\hat{E}_i, \; i=1,\ldots,m.$ There are two types of terms appearing on the RHS of (\ref{full}) and (\ref{reduced}): 
\begin{itemize}
\item 
i) terms that appear in both (\ref{full}) and (\ref{reduced}) with the same coefficient and 
\item
ii) terms that appear only on one of (\ref{full}) or (\ref{reduced}) or appear in both with different coefficients. 
\end{itemize} 
Terms of type i) should give rise to eigenvalues with value 1, for as long as the full system remains well resolved. As the full system starts to lose resolution (note that this is a gradual process due to the finiteness of the velocity of  propagation of disturbances across the Fourier modes), an eigenvalue due to a type i) term should start deviating from the value 1. Since all our calculations are performed in double precision, this statement can be made more accurate by looking to within how many digits of accuracy is the eigenvalue equal to 1.

Eigenvalues corresponding to terms of type ii) can have more varied behavior. However, there are some 
common aspects in the behavior of all these eigenvalues. The reason is that they come from terms that encode the differences between the full and the reduced system. 

\begin{itemize}
\item
During the time interval that the reduced set of modes $F$ contains all the energy, the differences between the full and reduced system should be zero, to within the precision of our simulations. As a result, matrix entries that come from terms that effect the drain of energy out of the resolved range should also be zero. It is easy to see, by looking at the definition of the matrix $B,$ that during the time interval when the reduced set of modes contains all the energy, the matrix $B$ will be (nearly) singular. Thus, the numerically computed values of the eigenvalues corresponding to these terms should just fluctuate wildly because the matrix $M$ is computed through the inversion of a (nearly) singular matrix.  
\item
After the initial time interval, energy starts seeping from $F$ to $G.$ This does not mean that the 
full system has become underresolved. The underresolution starts when the energy cannot escape from $G$ to modes corresponding to even smaller scales. During the interval when energy is moving from $F$ to $G,$ the conditions (\ref{conditions}) still hold but with a decreasing number of digits of accuracy. This means that both the reduced model and the reduced set of modes for the full system, lose energy at practically, but not exactly, the same rate. This is an important point . Clearly, the larger the number of modes available, the closer we can 
come to the singularity before energy starts to drain from $F$ to $G.$ 
\item
As time progresses, energy continues to flow from $F$ to $G$ bringing the full system closer and closer to underresolution. During that time the eigenvalues corresponding to type ii) terms increase. This increase signifies the deviation of the reduced model from the full system. If we could afford an infinite resolution, then the reduced and full systems would be the same until the singularity instant when their deviation would become infinite. Thus, the eigenvalue would stay zero until the singularity instant when it would become infinite. For any finite calculation, the instantaneous shift of the eigenvalue from zero to inifinity will be replaced by a period of rapid increase of the eigenvalue, followed by a turning point, and a period of slower increase of the eigenvalue until a finite maximum is reached. As we increase the resolution, the time of occurrence of the turning point should converge to the time of occurrence of the singularity and the maximum should grow unbounded.
\item
It is important to monitor to within how many digits of accuracy are conditions (\ref{conditions}) satisfied as we approach the eigenvalue turning point. As we have mentioned before, the correct digits of accuracy signify how reliable are the predictions for the temporal evolution of the eigenvalue. More details about this issue are provided in Section \ref{numerics} where we present specific numerical examples. 
\end{itemize}
Keeping in mind the subtleties mentioned above, we can formulate a simple algorithm for the determination of the occurrence of a possibly singular solution: 
\vskip14pt
{\bf Algorithm 1}
\begin{enumerate}
\item
Choose a resolution, i.e. a set of modes $F \cup G$ where $F$ are the resolved modes and $G$ the unresolved modes. Evolve simulatenously the full system for the modes in $F \cup G$ and 
a reduced model for the modes in $F.$
\item
Monitor the evolution of the eigenvalues of the renormalization matrix $M,$ as well as the digits of accuracy to within the conditions (\ref{conditions}) are satisfied. 
\item
The eigenvalues corresponding to type i) terms should be 1 until the reduced system starts deviating from the full system. After an initial phase when the matrix $B$ is singular (to within the numerical precision used), the eigenvalues corresponding to type ii) terms should increase until they reach a maximum. Record the time evolution of the eigenvalues up to time of the maximum. This will guarantee that the time of occurrence of the turning point is included.
\item 
Repeat the experiment with a higher resolution. If the solution of the PDE develops a singularity at time $T^{*},$ the sequence (as we increase the resolution) of the instants of occurrence of the turning point should converge to $T^{*}.$ The converse is not necessarily true. Convergence of the time of occurrence of the eigenvalue turning point to a finite value, say $\tau,$ does not imply that the solution 
of the PDE develops a singularity at $\tau.$ It may well be that the solution remains smooth but requires higher resolution than currently available. 
\end{enumerate} 
We repeat here that the accuracy of the predictions of the above algorithm relies on the accuracy of the reduced model used. We have chosen not to address this issue here, since our 
objective is to demonstrate how dimensionally reduced models can be used to study the behavior of solutions of PDEs and not how to obtain these dimensionally reduced models.


\subsection{How to approach efficiently a possible singularity}\label{approach}

The construction of the previous section suggests a way of performing mesh refinement. A progressive increase in the resolution, instead of 
starting with the largest possible resolution can allow us to reach the 
time when the interesting dynamics occur more efficiently. In particular, the entries of matrix $B$ provide a way of monitoring the transfer of energy from $F$ to $G.$ When all the energy is contained in $F,$ the matrix $B$ is singular (to within the numerical precision used). When energy starts transferring from $F$ to $G,$ the absolute value of the determinant of $B,$ $|detB|,$ acquires positive values. One can think of $|detB|$ as an error control criterion, which can be used to decide when it is time to refine the mesh. Since we are restricting ourselves to the case of periodic boundary conditions and Fourier modes, the mesh refinement is performed in Fourier space, by appending to the vector of Fourier modes, a vector of higher Fourier modes with zero amplitude (this procedure is known as spectral interpolation \cite{boyd}). For example, when $|detB| \geq TOL,$ where $TOL$ a specified tolerance we can double the number of Fourier modes in each spatial direction. This leads to a uniform refinement of the mesh in real space. A real space formulation of the algorithms which leads to a real space formulation of the mesh refinement criterion and allows the treatment of non-periodic boundary conditions and more complicated geometries will be presented elsewhere.

We discuss now a way to calibrate the mesh refinement criterion $TOL.$ The strictest value allowed is to set $TOL=\eps,$ where $\eps$ the numerical precision used e.g. double precision. This means that we refine the moment we detect the smallest transfer of energy (allowed by our numerical precision) from $F$ to $G.$ This can be unnecessarily costly because at any given moment we do not allow any energy to be present in the modes in $G.$ However, we know (see also discussion in the previous section), that the modes in $G$ are allowed to have some energy {\it without} our calculations becoming underresolved. Thus, $TOL$ can be allowed to have values larger than $\eps$, and still get a mesh refinement scheme that does not introduce any substantial error. One way to calibrate the value of 
$TOL$ is the following:
\vskip14pt
{\bf Algorithm 2}
\begin{enumerate}
\item
Choose a value for $TOL.$ For this value of $TOL$ run a mesh refinement calculation, starting, say, from $N_{start}$ modes to $N_{final}$ modes. For example, let $N_{start}=32$ and double at each refinement until, say $N_{final}=256$ modes. Record the values of the quantities $\hat{E}_i, \; i=1,\ldots,m$ when $N=N_{final}$ and $|detB|=TOL.$ Let's call this simulation $S1.$
\item
For the same value of $TOL$ run a calculation with $N_{start}=N_{final}$ modes (for the example $N_{start}=N_{final}=256$). Record the values of the quantities $\hat{E}_i, \; i=1,\ldots,m$ when $|detB|=TOL.$ Let's call this simulation $S2.$
\item
Compare to within how many digits of accuracy the quantities $\hat{E}_i, \; i=1,\ldots,m$ computed from $S1$ and $S2$ agree. If the agreement is to within a specified accuracy, say 5 digits, then choose 
this value of $TOL.$ If the agreement is in fewer digits, then decrease $TOL$ (more stringent criterion) and repeat until agreement is met.
\item
Use the above decided value of $TOL$ to perform a mesh refinement calculation with a larger magnification ratio, i.e. a larger value for the ratio $N_{final}/ N_{start}.$  
\end{enumerate}
We want to make two comments about the proposed scheme. First, there is no a priori guarantee that the use of the same value of $TOL$ for a mesh refinement calculation with a larger magnification ratio, will preserve the number of digits of accuracy. However, the fact that we choose the value of $TOL$ based on the rate of energy transfer from $F$ to $G$ and {\it not} on how much energy is contained in $G$ makes us  hopeful that the proposed refinement scheme has a sound basis. Second, the scheme relies on the accuracy of the reduced model since it relies on the determinant of $B.$ As before, we assume that we have access to a good reduced model that will detect correctly the transfer of energy from $F$ to $G.$


\subsection{How to follow a possible singularity}\label{follow}
We conclude our discussion of the algorithms with an algorithm that allows the tracking of a possible singularity:
\vskip14pt
{\bf Algorithm 3}
\begin{enumerate}
\item
Use {\bf Algorithm 2} until the maximum allowed resolution has been reached. Suppose that the set of modes for the maximum resolution 
is $K.$
\item
Divide $K$ in two sets $F$ and $G,$ where $F$ the resolved variables and $G$ the unresolved variables. Use the reduced model for the modes in $F$ to continue the calculation to later times.
\end{enumerate}
We should note, that from the numerical analysis perspective, we are interested in following correctly a solution whether it is singular or only near-singular. Even if future advancements in computational power and/or mathematical analysis allow us to decide on the existence or not of a singularity, it is still important to have a good reduced model which can reproduce, at a lower cost, the essential features of the solution.



\subsection{Numerical results for Case I}\label{numerics}
We present numerical results of the three algorithms discussed above for the 1D Burgers equation with zero and nonzero viscosity and periodic boundary conditions. Our purpose in this paper is to demonstrate the use of the algorithms. Extension of the algorithms to non-periodic boundary conditions, more spatial dimensions, systems of partial differential equations and more complicated geometries will be presented elsewhere. 

The 1D Burgers equation is given by 
\begin{equation}\label{burgers}
u_t+u u_x = \nu u_{xx}
\end{equation}
where $\nu$ is the viscosity coefficient. Equation (\ref{burgers}) should be supplemented with an initial condition $u(x,0)=u_0(x)$ and boundary conditions. We solve (\ref{burgers}) in the interval $[0,2\pi]$ with periodic boundary conditions. This allows us to expand the solution in Fourier series
$$v_{M}(x,t )=\underset{k \in F \cup G}{\sum} u_k(t) e^{ikx},$$
where $F \cup G=[-\frac{M}{2},\frac{M}{2}-1].$ We have written the set of Fourier modes as the union of two sets 
in anticipation of the construction of the reduced model comprising only of the modes in $F=[-\frac{N}{2},\frac{N}{2}-1],$ where $ N < M.$
The equation of motion for the Fourier mode $u_k$ becomes
\begin{equation}
\label{burgersode}
 \frac{d u_k}{dt}=- \frac{ik}{2} \underset{p, q \in F \cup G}{\underset{p+q=k  }{ \sum}} u_{p} u_{q} - \nu k^2 u_k.
\end{equation}

\subsubsection{The $t$-model}\label{numericstmodel}

We need to choose a reduced model for the modes in $F.$ We use a reduced model, known as the $t$-model, which follows correctly the behavior of the solution to the inviscid Burgers equation even after the formation of shocks \cite{bernstein,HS06}. The $t$-model was first derived in the context of statistical irreversible mechanics \cite{CHK3} and was later analyzed in \cite{bernstein,HS06}. It is based on the assumption of the absence of time scale separation between the resolved and unresolved modes. We will use the same model for the case with nonzero viscosity and comment on its validity when appropriate. For a mode $u_k'$ in $F$ the model is given by
\begin{multline}\label{burgersode2}
\frac{d}{dt}u_k'=- \frac{ik}{2}   \underset{p \in F ,\, q \in F }{\underset{p+q=k  }{ \sum}}   u_{p}'u_{q}'  -\nu k^2 u_k'\\
-\frac{ik}{2} \underset{p \in F  ,\, q \in G}{\underset{p+q=k  }{ \sum}}   u_{p}'  \biggl[ -t\frac{iq}{2}  \underset{r \in F ,\,  s \in F }{\underset{r+s=q }{ \sum}}  u_{r}'u_{s}' \biggr] \\
-\frac{ik}{2} \underset{p \in G ,\, q \in F }{\underset{p+q=k  }{ \sum}}  \biggl[ -t\frac{ip}{2} \underset{r \in F  ,\, s \in F }{\underset{r+s=p }{ \sum}} u_{r}'
u_{s}' \biggr]  u_{q}' . 
\end{multline}
The first term on the RHS of (\ref{burgersode2}) is of the same form as the first term in (\ref{burgersode}), except that the term in (\ref{burgersode2}) is defined only for the modes in $F.$ The viscous term is the same. The third and fourth terms in (\ref{burgersode2}) are not present in (\ref{burgersode}). They are cubic in the Fourier modes and they are effecting the drain of energy out of the modes in $F.$ We should note here that the cubic terms in the $t$-model do not depend on the viscosity. To conform with the notation in Section \ref{algo} we rewrite (\ref{burgersode2}) as
\begin{multline*}
\frac{d}{dt}u_k'= a^{(1)}_1 \biggl[- \frac{ik}{2}   \underset{p \in F ,\, q \in F }{\underset{p+q=k  }{ \sum}}   u_{p}'u_{q}'  -\nu k^2 u_k' \biggr] +\\
a^{(1)}_2 \Biggl[ -\frac{ik}{2} \underset{p \in F  ,\, q \in G}{\underset{p+q=k  }{ \sum}}   u_{p}'  \biggl[ -t\frac{iq}{2}  \underset{r \in F ,\,  s \in F }{\underset{r+s=q }{ \sum}}  u_{r}'u_{s}' \biggr] \\
-\frac{ik}{2} \underset{p \in G ,\, q \in F }{\underset{p+q=k  }{ \sum}}  \biggl[ -t\frac{ip}{2} \underset{r \in F  ,\, s \in F }{\underset{r+s=p }{ \sum}} u_{r}'
u_{s}' \biggr]  u_{q}' \Biggr],  
\end{multline*}
where $a^{(1)}_1=1$ and $a^{(1)}_2=1.$ We rewrite Equation (\ref{burgersode})  as
\begin{multline*}
 \frac{d u_k}{dt}= a^{(0)}_1 \biggl[ - \frac{ik}{2} \underset{p, q \in F \cup G}{\underset{p+q=k  }{ \sum}} u_{p} u_{q} - \nu k^2 u_k \biggr] + \\
a^{(0)}_2 \Biggl[ -\frac{ik}{2} \underset{p \in F\cup G  ,\, q \in I}{\underset{p+q=k  }{ \sum}}   u_{p}  \biggl[ -t\frac{iq}{2}  \underset{r \in F\cup G ,\,  s \in F \cup G }{\underset{r+s=q }{ \sum}}  u_{r}u_{s} \biggr] \\
-\frac{ik}{2} \underset{p \in I ,\, q \in F\cup G }{\underset{p+q=k  }{ \sum}}  \biggl[ -t\frac{ip}{2} \underset{r \in F\cup G  ,\, s \in F\cup G }{\underset{r+s=p }{ \sum}} u_{r} u_{s} \biggr]  u_{q} \Biggr]  ,
\end{multline*}
where $a^{(0)}_1=1$ and $a^{(0)}_2=0.$ The reader should note that we have introduced a new set of modes $I.$ This is the set of unresolved modes for the {\it full} system. The reason for introducing the set $I$ is that, as is the case in renormalization formulations, the terms appearing in the RHS of the equations at the different levels of resolution should be of the same functional form. The difference between the different levels of resolution should be only in the range of modes used. Since the $t$-model involves a quadratic convolution sum with one index in the resolved range and the other in the unresolved range, we should use the same functional form when constructing the corresponding term for the full system. Thus, this term should involve a convolution sum with one index in the range $F \cup G$ and the other in $I.$ 

Further, define
$$\hR^{(0)}_{1k}(t,\hat{u}(t))=- \frac{ik}{2} \underset{p, q \in F \cup G}{\underset{p+q=k  }{ \sum}} u_{p} u_{q} - \nu k^2 u_k $$
and
\begin{multline*}
\hR^{(0)}_{2k}(t,\hat{u}(t))= -\frac{ik}{2} \underset{p \in F\cup G  ,\, q \in I}{\underset{p+q=k  }{ \sum}}   u_{p}  \biggl[ -t\frac{iq}{2}  \underset{r \in F\cup G ,\,  s \in F \cup G }{\underset{r+s=q }{ \sum}}  u_{r}u_{s} \biggr] \\
-\frac{ik}{2} \underset{p \in I ,\, q \in F\cup G }{\underset{p+q=k  }{ \sum}}  \biggl[ -t\frac{ip}{2} \underset{r \in F\cup G  ,\, s \in F\cup G }{\underset{r+s=p }{ \sum}} u_{r} u_{s} \biggr]  u_{q}
\end{multline*}
Also, define
$$\hR^{(1)}_{1k}(t,\hat{u}'(t))=- \frac{ik}{2} \underset{p, q \in F }
{\underset{p+q=k  }{ \sum}} u_{p}' u_{q}' - \nu k^2 u_k' $$
and
\begin{multline*}
\hR^{(1)}_{2k}(t,\hat{u}'(t))=-\frac{ik}{2} \underset{p \in F  ,\, q \in G}{\underset{p+q=k  }{ \sum}}   u_{p}'  \biggl[ -t\frac{iq}{2}  \underset{r \in F ,\,  s \in F }{\underset{r+s=q }{ \sum}}  u_{r}'u_{s}' \biggr] \\
-\frac{ik}{2} \underset{p \in G ,\, q \in F }{\underset{p+q=k  }{ \sum}}  \biggl[ -t\frac{ip}{2} \underset{r \in F  ,\, s \in F }{\underset{r+s=p }{ \sum}} u_{r}'
u_{s}' \biggr]  u_{q}'
\end{multline*}
Thus, the equations of motion for the resolved modes in the full system and the reduced model can be written as
\begin{equation}\label{alburgersode}
 \frac{d u_k}{dt}=\sum_{i=1}^{2} a^{(0)}_i \hR^{(0)}_{ik}(t,\hat{u}(t))
\end{equation}
and
\begin{equation}\label{alburgersode2}
\frac{d u_k'}{dt}=\sum_{i=1}^{2} a^{(1)}_i \hR^{(1)}_{ik}(t,\hat{u}'(t))
\end{equation}
We want to make two comments about the rewriting of Equations (\ref{burgersode}) and (\ref{burgersode2}). First, we have $a^{(0)}_1=a^{(1)}_1=1.$ This term should give rise to an eigenvalue equal to 1 for the matrix $M$ as explained in Section \ref{location}.   Second, we have $a^{(1)}_2 \neq a^{(0)}_2,$ so the corresponding eigenvalue should be different than 1.  

To proceed with the algorithms we need to define the quantities $\hat{E}_i, \; i=1,\ldots,m.$ In our case, $m=2$ and we need to define $\hat{E}_1$ and $\hat{E}_2.$ The choice of the $\hat{E}_i$ is not unique. We chose for our experiments $\hat{E}_1=\underset{k \in F}{\sum} |u_k|^2$ and $\hat{E}_2=\underset{k \in F}{\sum} |u_k|^4.$ The rates of change of the $\hat{E}_i$ are given for the full system by
$$\frac{d\hat{E}_1}{dt}=\underset{k \in F}{\sum}a^{(0)}_1 2 Re(\hR^{(0)}_{1k}(t,\hat{u}(t)) u_k^{*})+ a^{(0)}_2 2 Re(\hR^{(0)}_{2k}(t,\hat{u}(t)) u_k^{*})$$
and
$$\frac{d\hat{E}_2}{dt}=\underset{k \in F}{\sum}a^{(0)}_1 2 Re(2 \hR^{(0)}_{1k}(t,\hat{u}(t)) |u_k|^2 u_k^{*})+ a^{(0)}_2 2 Re(2 \hR^{(0)}_{2k}(t,\hat{u}(t)) |u_k|^2 u_k^{*})$$
where 
$u_k^{*}$ is the complex conjugate of $u_k.$ Similarly, for the reduced system we have
$$\frac{d\hat{E}_1}{dt}=\underset{k \in F}{\sum}a^{(1)}_1 2 Re(\hR^{(1)}_{1k}(t,\hat{u}'(t)) u_k'^{*})+ a^{(1)}_2 2 Re(\hR^{(1)}_{2k}(t,\hat{u}'(t)) u_k'^{*})$$
and
$$\frac{d\hat{E}_2}{dt}=\underset{k \in F}{\sum}a^{(1)}_1 2 Re(2 \hR^{(1)}_{1k}(t,\hat{u}'(t)) |u_k'|^2 u_k'^{*})+ a^{(1)}_2 2 Re(2 \hR^{(1)}_{2k}(t,\hat{u}'(t)) |u_k'|^2 u_k'^{*})$$ 
The equations for the rates of change of the $\hat{E}_i$ can be used for the computation of the $2 \times 2$ matrices $A$ and $B$ through the relations (\ref{rngmatrix1}) of Section \ref{algo}. The matrix $M$ can be computed through the matrix equation (\ref{rngmatrix2}).

\paragraph{The case of zero viscosity}\label{zero}

We begin our study of the algorithms with the case of zero viscosity, i.e. $\nu=0.$ We choose the initial condition $u(x,0)=u_0(x)=\sin x,$ which develops a shock at $x=\pi / 2$ at $T^{*}=1.$ No matter how large a resolution we use, a simulation based on (\ref{burgersode}) will become underresolved before $T^{*}=1.$ Of course, the larger the resolution the closer the simulation will get to the time of the shock formation before it becomes underresolved. The solution before the formation of the shock conserves its $L_2$ norm (energy). After the formation of the shock, the energy of the solution decays. A simulation based on (\ref{burgersode}) conserves energy for all time. On the other hand, a simulation based on (\ref{burgersode2}) remains well resolved and dissipates energy at the correct rate \cite{bernstein,HS06}.

\begin{figure}
\centering
\subfigure[]{\epsfig{file=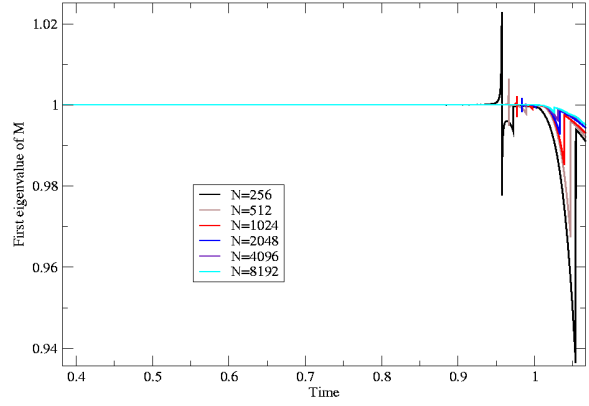, width=3.in}}
\qquad
\vskip14pt
\subfigure[]{\epsfig{file=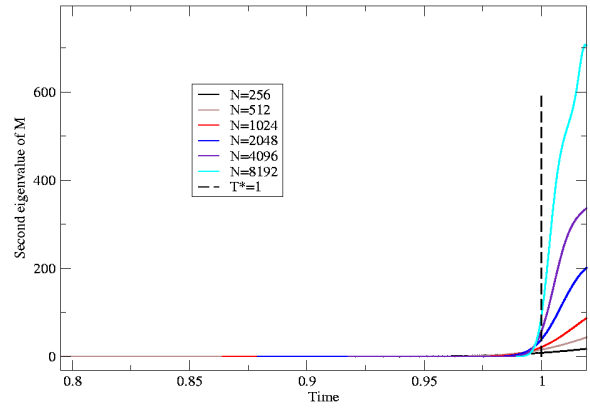,width=3.in}}
\caption{Case I, $\nu=0$: (a) Evolution of the first eigenvalue of the matrix $M$ for different resolutions. (b) Evolution of the second eigenvalue.}
\label{plot_a1v0}
\end{figure}

\begin{figure}
\centering
\epsfig{file=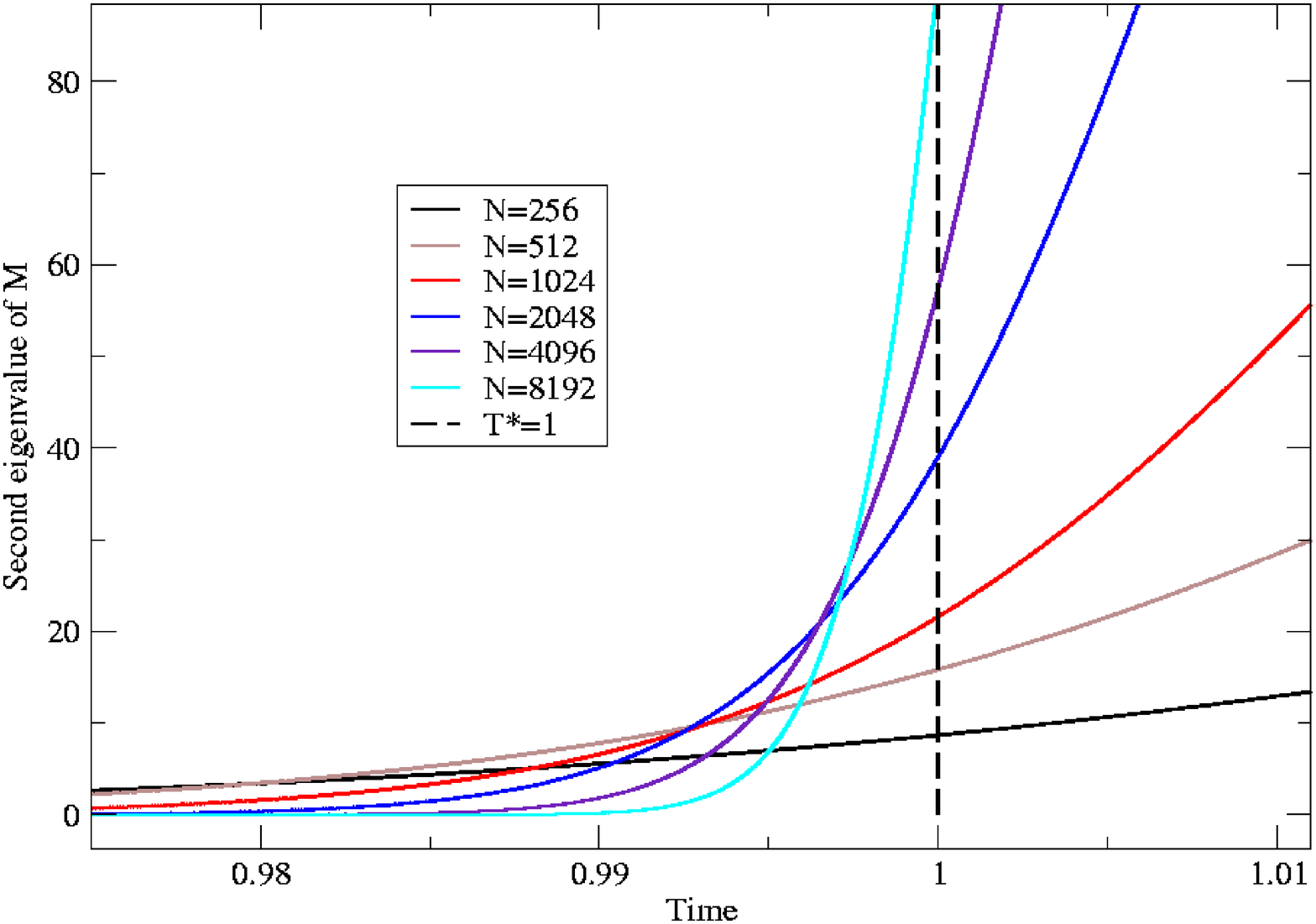,width=3.in}
\caption{Case I, $\nu=0$: Detail of the evolution of the second eigenvalue around the time of singularity. As the resolution is increased, the time when the evolution with a higher
resolution crosses the evolution with the lower resolution, moves towards 1.}
\label{plot_a1v0_detail_a}
\end{figure}

\begin{figure}
\centering
\epsfig{file=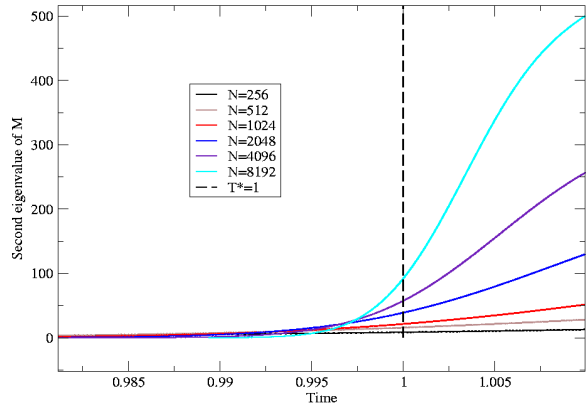,width=3.in}
\caption{Case I, $\nu=0$: Detail of the evolution of the second eigenvalue around the time of singularity. The turning point moves towards 1 as the resolution is increased.}
\label{plot_a1v0_detail_b}
\end{figure}

First, we present results of the application of Algorithm 1 (from now one referred to as A1). We use this algorithm to locate the time of occurrence of a singularity. We repeat that for this example we know that a singularity occurs at $T^*=1$ and thus, the algorithm should show that. If we did not know that a singularity exists, the algorithm should suggest that something is happening but we would not know whether it is a singular solution or a near singular solution that just needs more modes to be resolved than we can afford. Figure \ref{plot_a1v0} shows the evolution of the first and second eigenvalues of the matrix $M$ for different resolutions. As we expected, the first eigenvalue is equal to 1 until the moment when the reduced and full system start deviating. The second eigenvalue remains close to zero until the time of deviation of the reduced and full system. The evolution of the second eigenvalue is more informative about the behavior of the two systems. As we increase the resolution, the second eigenvalue remains closer to zero for a longer time. When it increases it does so at a rate that is higher the higher is the resolution. In Figures \ref{plot_a1v0_detail_a} and \ref{plot_a1v0_detail_b}  we have zoomed in on the time interval [.985,1.005]. As it should happen, the time when the evolution of the second eigenvalue for a higher resolution crosses the evolution of the second eigenvalue for the preceding resolution, moves towards 1. Also, the time of occurrence of the turning point moves towards 1 as we increase the resolution. The turning point for the highest resolution here $N=8192$ is around 1.003. As we have already explained, one should keep track of the accuracy with which the conditions \eqref{conditions} are satisfied. For $N=256,$ the conditions are satisfied to 5 significant digits at time $T^{*}=1.$ For $N=8192,$ the accuracy has increased to 8 significant digits.

\begin{figure}
\centering
\subfigure[]{\epsfig{file=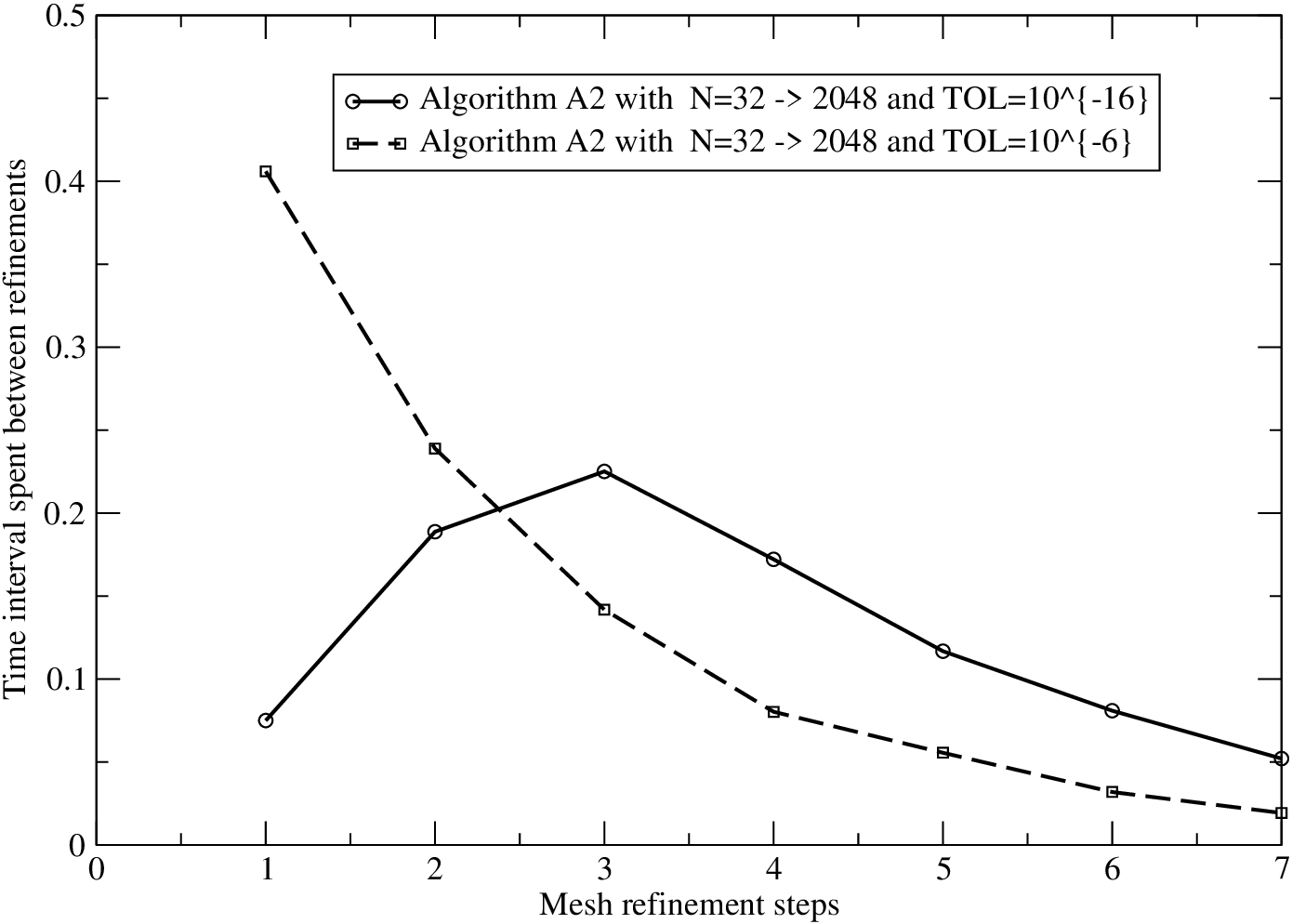,width=3.in}}
\qquad
\vskip20pt
\subfigure[]{\epsfig{file=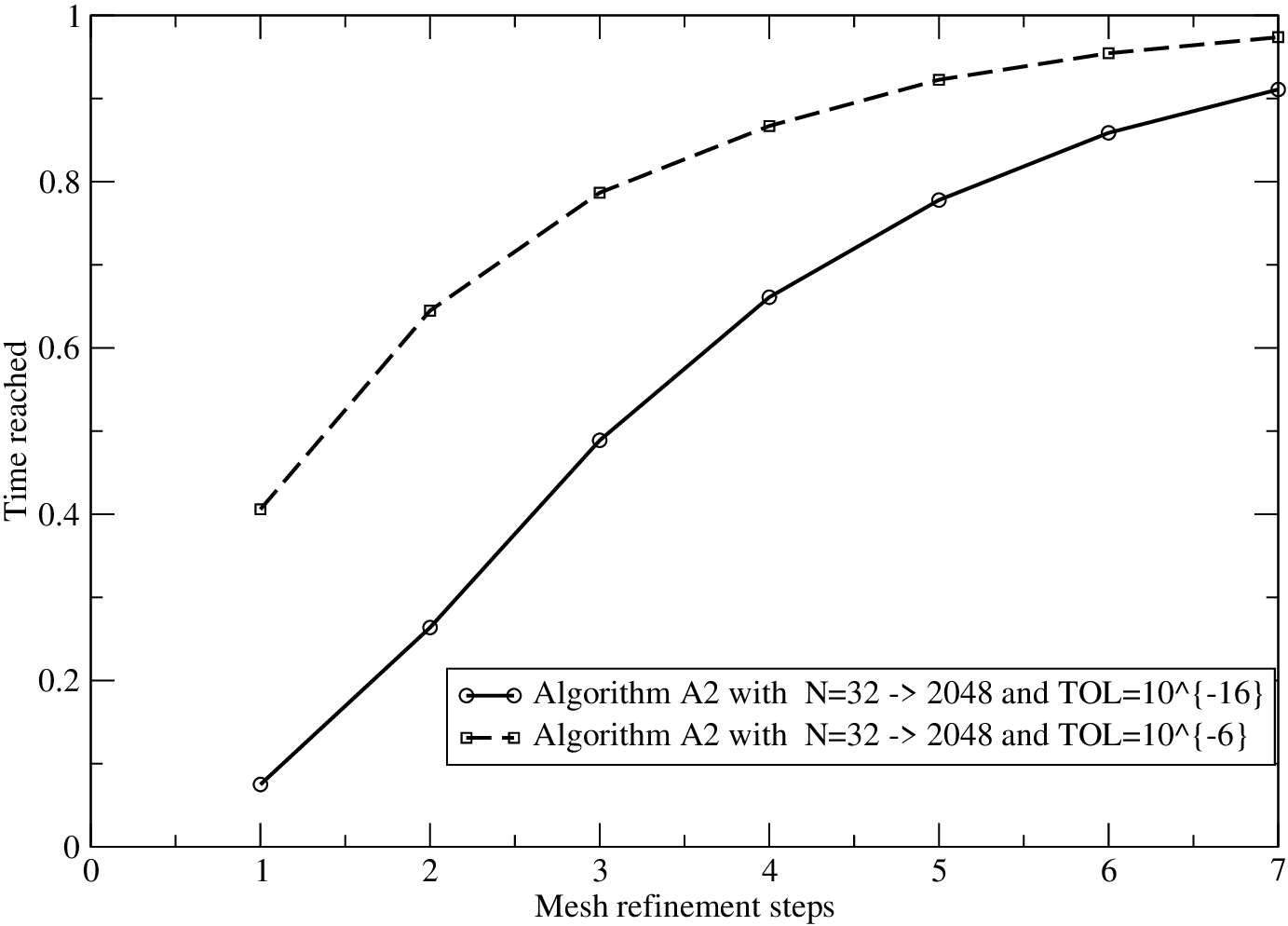,width=3.in}}
\caption{Case I, $\nu=0$: (a) Time spent between refinement steps for different tolerance values. (b) Time reached with the maximum allowed resolution. }
\label{plot_a2v0}
\end{figure}

The results of the Algorithm 2 (A2) are presented in Figure \ref{plot_a2v0}. We use this algorithm to perform mesh refinement when needed starting from a resolution $N_{start}=32$ and allowing a maximum resolution of $N_{final}=2048.$ We present results for two values of the tolerance $TOL1=10^{-16}$ and $TOL2=10^{-6}.$ For the case of $TOL1,$ where we set the tolerance equal to double precision, the first  few refinement steps show an increase of the time spent between refinements. However, in both cases, the time spent between successive refinements eventually starts shrinking with the number of refinement steps. We also plot the time reached with the maximum allowed resolution. When the tolerance criterion is less strict the algorithm can reach later times before running out of resolution. 

\begin{figure}
\centering
\epsfig{file=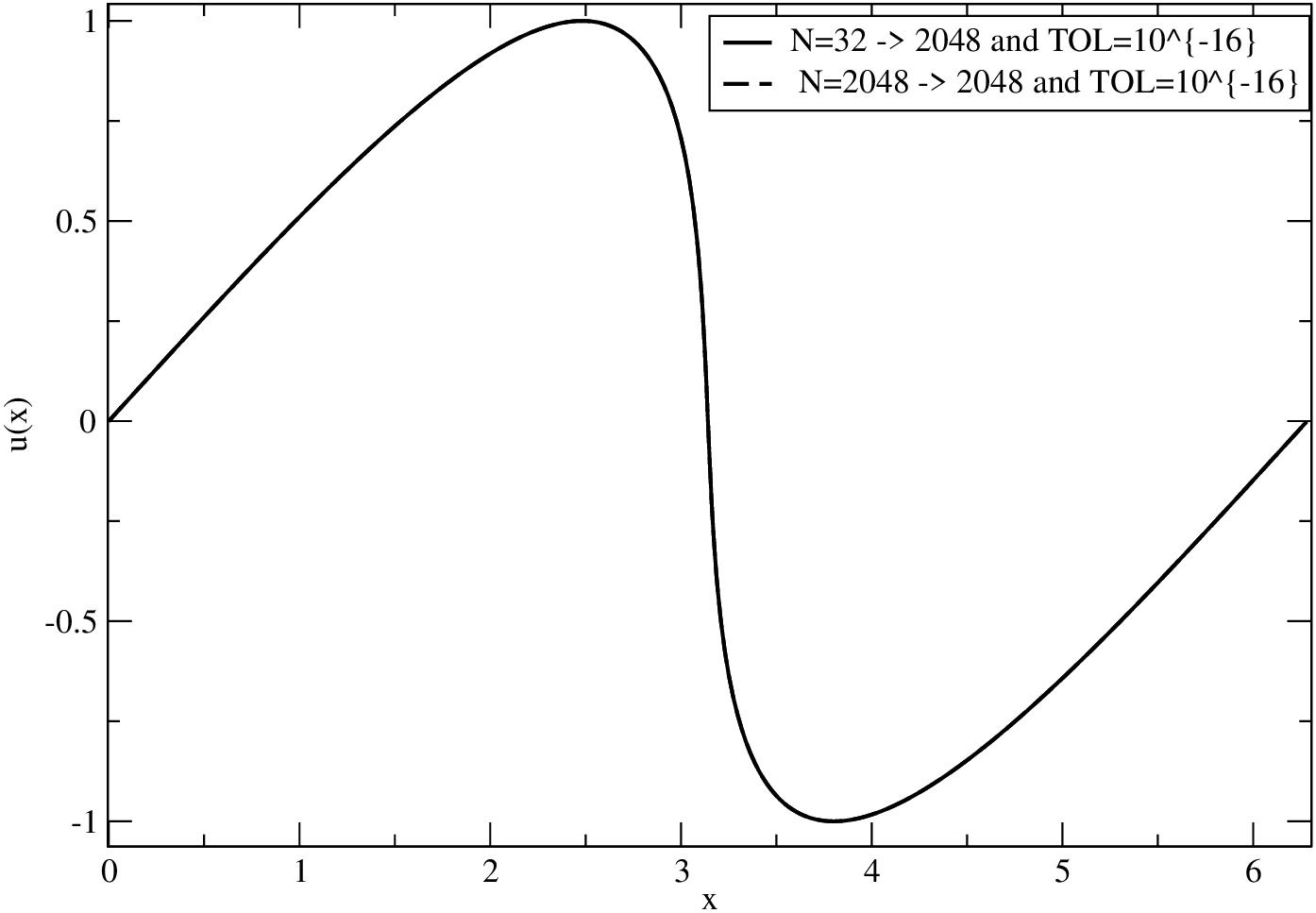,width=3.in}
\caption{Case I, $\nu=0$: Comparison of the velocity field produced at the time of termination of A2 for two different magnification ratios. 
The first simulation has $N_{start}=32$ and $N_{final}=2048$ while the second has $N_{start}=N_{final}=2048.$ }
\label{plot_a2v0_velocity}
\end{figure}

In Figure \ref{plot_a2v0_velocity} we compare the velocity field produced by A2 with $N_{start}=32,$ $N_{final}=2048$ and $TOL1=10^{-16}$ with the velocity field produced by A2 with $N_{start}=N_{final}=2048$ and the same tolerance. It is obvious that the results are in very good agreement. However, the mesh refinement calculation was about 9 times faster to run. The acceleration factor, with practically the same accuracy, can be increased by using a larger magnification ratio.

Next, we present results for the Algorithm 3 (A3). The algorithm uses A2 until the maximum resolution is reached and then switches to the t-model. We compare the energy evolution as predicted by A3 with $N_{start}=32, \; N_{final}=512$ and $TOL=10^{-16}$, a $t$-model calculation with $N=512$ and the random choice method \cite{chorinran} with $N=4096$ points. The mesh refinement algorithm produces results that are in very good agreement with the correct energy decay as computed by the random choice method.    

\begin{figure}
\centering
{\epsfig{file=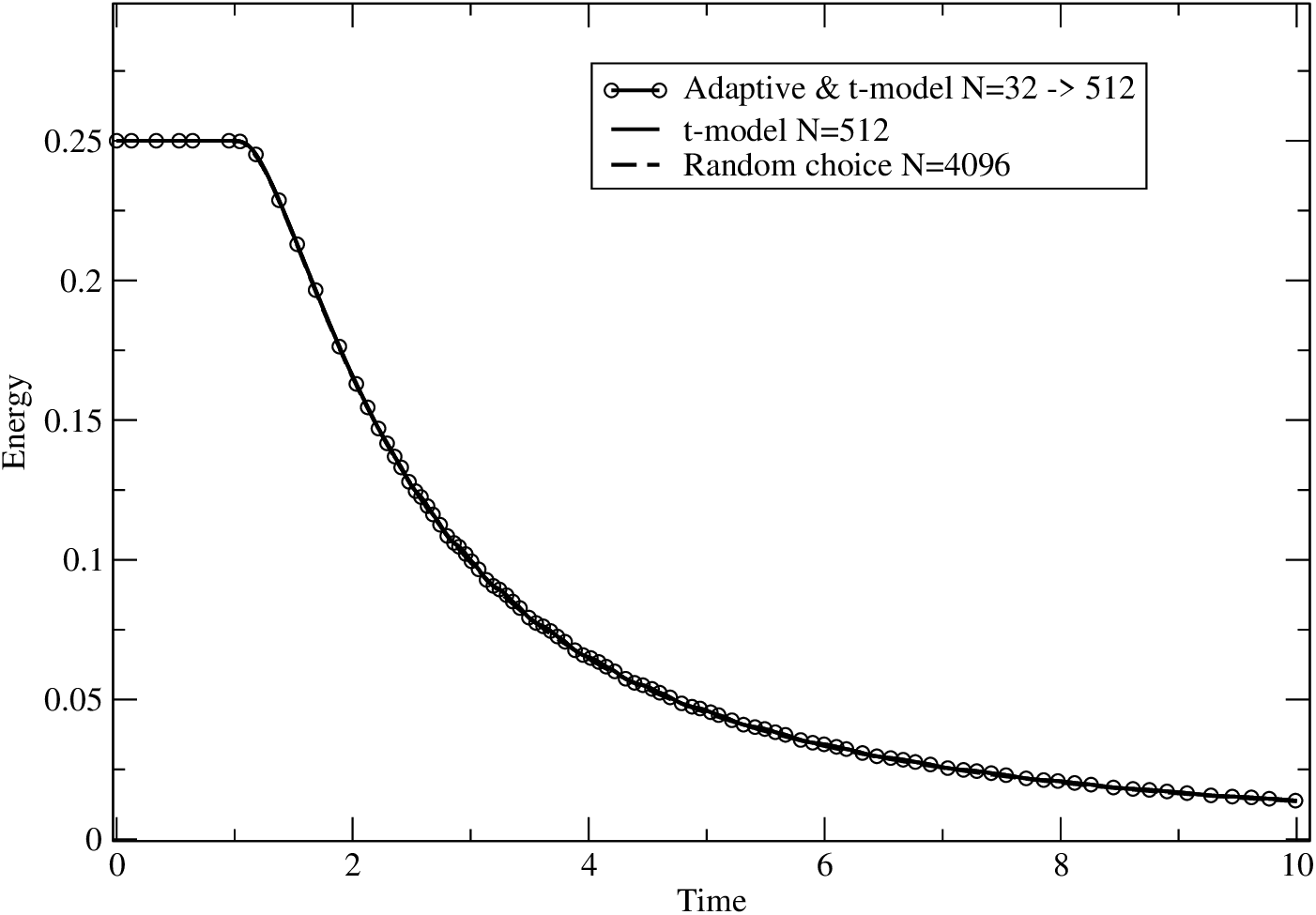,width=3.in}}
\caption{Case I, $\nu=0$: Comparison of the energy evolution as produced by A3, the $t$-model and the random choice method.}
\label{plot_a3v0_energy}
\end{figure}

\paragraph{The case of nonzero viscosity}\label{nonzero}

\begin{figure}
\centering
{\epsfig{file=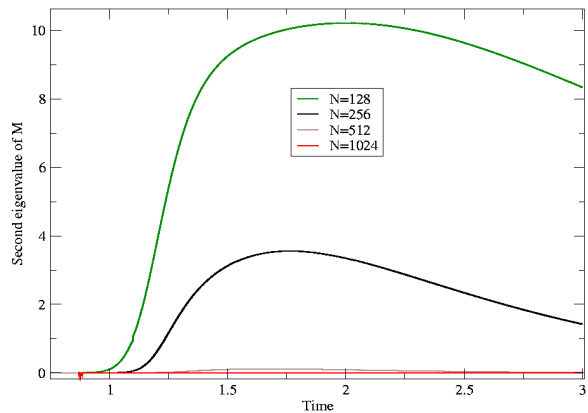,width=3.in}}
\caption{Case I, $\nu=0.01$: Evolution of the second eigenvalue of the matrix $M$ for different resolutions.}
\label{plot_a1v_2}
\end{figure}

\begin{figure}
\centering
\subfigure[] {\epsfig{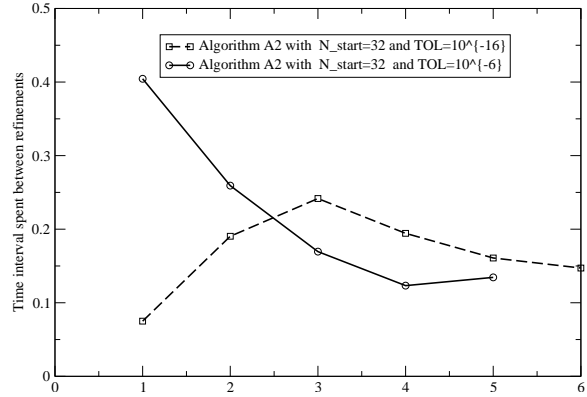}}
\qquad
\subfigure[] {\epsfig{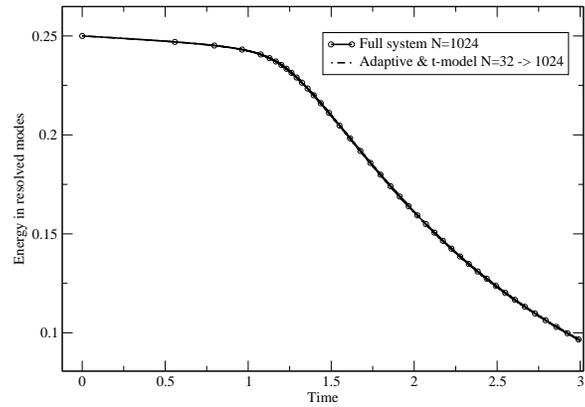}}
\caption{Case I, $\nu=0.01$: (a) Time spent between refinement steps for different tolerance values. (b) Comparison of the energy evolution as produced by A3 and the full system.}
\label{plot_a2v_2}
\end{figure}

When the viscosity is nonzero the solution remains smooth for all time. However, depending on the value of the viscosity there is a certain number of modes needed to fully resolve the solution. Figure \ref{plot_a1v_2} shows the evolution of the second eigenvalue of the matrix $M$ for different resolutions with $\nu=0.01.$ Notice that the behavior of the eigenvalue is different from the case of no viscosity. There the eigenvalue increased with an increase in the resolution. Here it is decreasing until for large enough resolution it is practically zero for the duration of the simulation. In Figure \ref{plot_a2v_2}(a) we present the progress of the refinement steps for $TOL1=10^{-16}$ and $TOL2=10^{-6}.$  We present the number of refinement steps needed until the refinement criterion is no longer activated. In the first case, the criterion is very stringent and the algorithm has to refine more times to achieve it. Figure \ref{plot_a2v_2}(b) shows the evolution of the energy as predicted by A3 with $N_{start}=32, \; N_{final}=1024$ and $TOL=10^{-16}$ and a full system with $N=1024$ modes. The algorithm A3 produces results that are in very good agreement with the full system and does so 3 times faster. 

\subsubsection{The Galerkin model}\label{numericsgalerkin}
We also present results for the case when the reduced model is the Galerkin model, i.e. where we set the unresolved modes equal to zero for all times. This means that $a^{(1)}_1=1, \; a^{(1)}_2=0.$ Note that in the limit of infinite resolution the Galerkin model tends to the Galerkin truncation for the full system. The algorithm A1 which computes the eigenvalues should be able to detect that. Since $a^{(0)}_2=a^{(1)}_2=0,$ we can restrict ourselves to one quantity $\hat{E}_1.$ However, if we choose the $L_2$ norm of the resolved modes as $\hat{E}_1$ will lead to numerical trouble, because  $\frac{d\hat{E}_1}{dt}=0$ for the reduced Galerkin model. That would lead to $|detB|=0$ to within the numerical precision used. We can remedy the numerical issues by using two quantities $\hat{E}_1$ and $\hat{E}_2$ as before. The difference with the case of the $t$-model is that here both the full and reduced models do {\it not} involve the cubic term (since $a^{(0)}_2=a^{(1)}_2=0$). Thus, the behavior of the second eigenvalue will be different. In fact, since we know that the reduced Galerkin model converges to the full Galerkin system, the second eigenvalue should decrease as we increase the resolution.

\begin{figure}
\centering
\subfigure[]{\epsfig{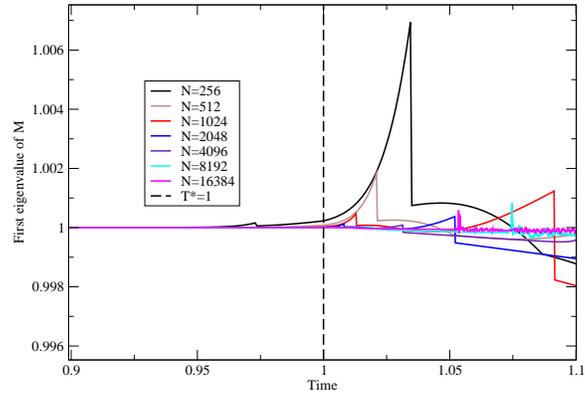}}
\qquad
\vskip14pt
\subfigure[]{\epsfig{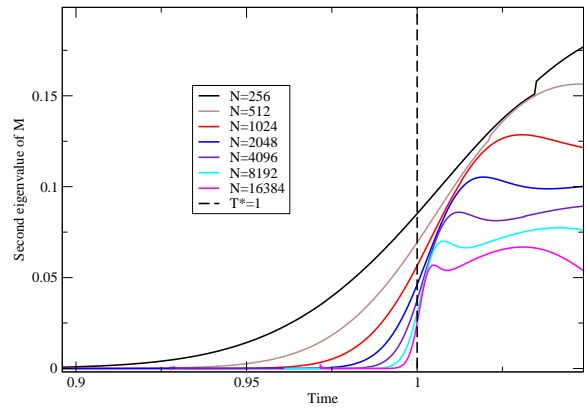}}
\caption{Case I, Galerkin, $\nu=0$: (a) Evolution of the first eigenvalue of the matrix $M$ for different resolutions. (b) Evolution of the second eigenvalue.}
\label{plot_a1v0_galerkin}
\end{figure}

\paragraph{The case of zero viscosity}

\begin{figure}
\centering
{\epsfig{file=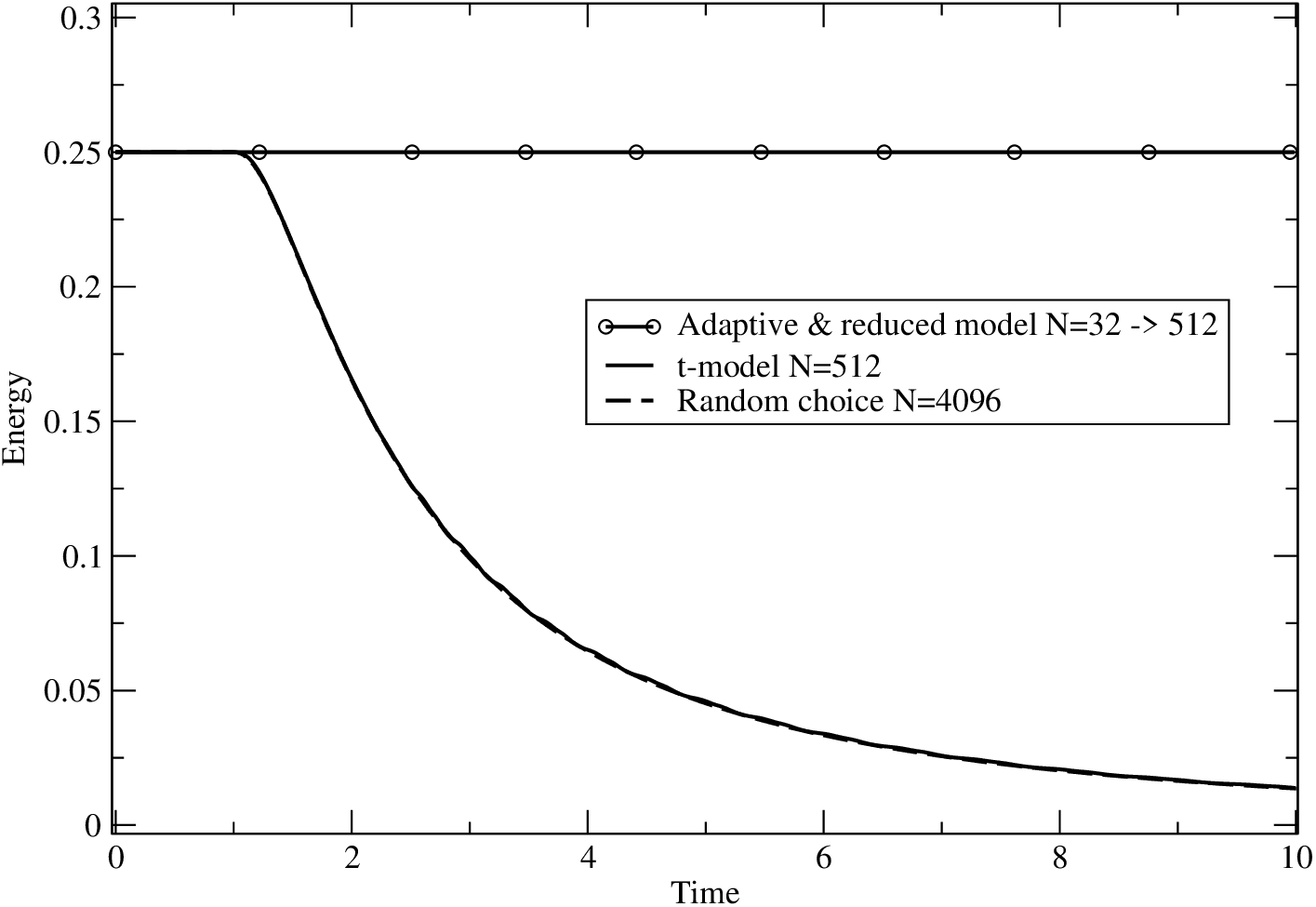,width=3.in}}
\caption{Case I, Galerkin, $\nu=0$: Comparison of the energy evolution as produced by A3, the $t$-model and the random choice method.}
\label{plot_a3v0_energy_galerkin}
\end{figure}

Figure \ref{plot_a1v0_galerkin} shows the evolution of the first and second eigenvalues of the matrix $M$ for different resolutions with $\nu=0.$ As we expected, the first eigenvalue is equal to 1 until the moment when the reduced and full system start deviating. As we increase the resolution, the first eigenvalue remains 1 for a longer time. When it increases it does so at a rate that is slower the higher is the resolution. Also, the second eigenvalue decreases as the resolution is increased. This means that the reduced model converges to the full system as the resolution is increased. For $N=256,$ the conditions \eqref{conditions} are satisfied to 4 significant digits at time $T^{*}=1.$ For $N=8192,$ the accuracy has increased to 7 significant digits. Note that compared to the case when the reduced model is the $t$-model the reduced model here loses one additional significant digit by the time of the singularity occurrence. Of course, as will be shown later (see Figure \ref{plot_a3v0_energy_galerkin}) the Galerkin reduced model is unable to dissipate energy and it quickly becomes very bad as a model for the weak solution (shock).  

\begin{figure}
\centering
\subfigure[]{\epsfig{file=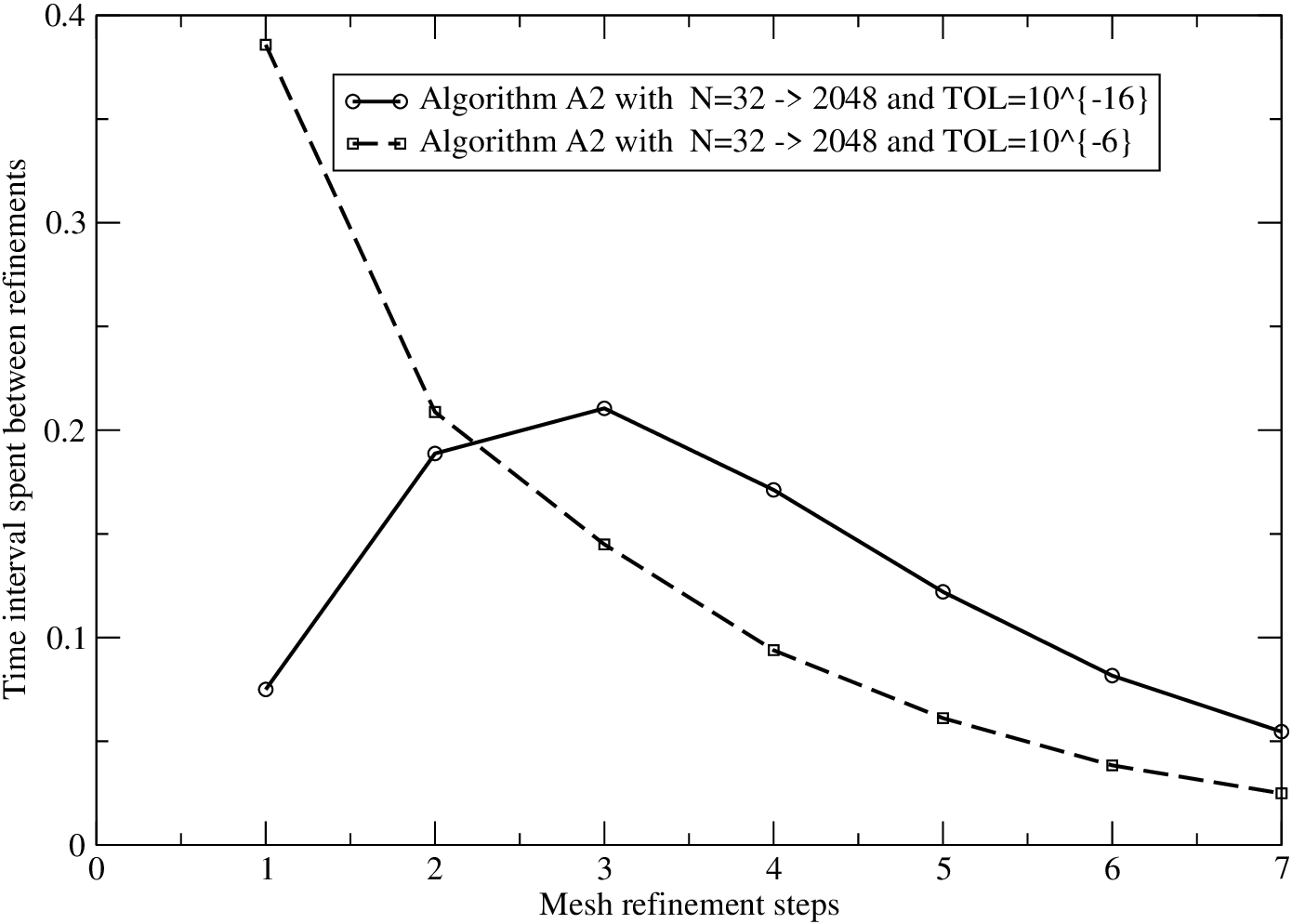,width=3.in}}
\qquad
\vskip20pt
\subfigure[]{\epsfig{file=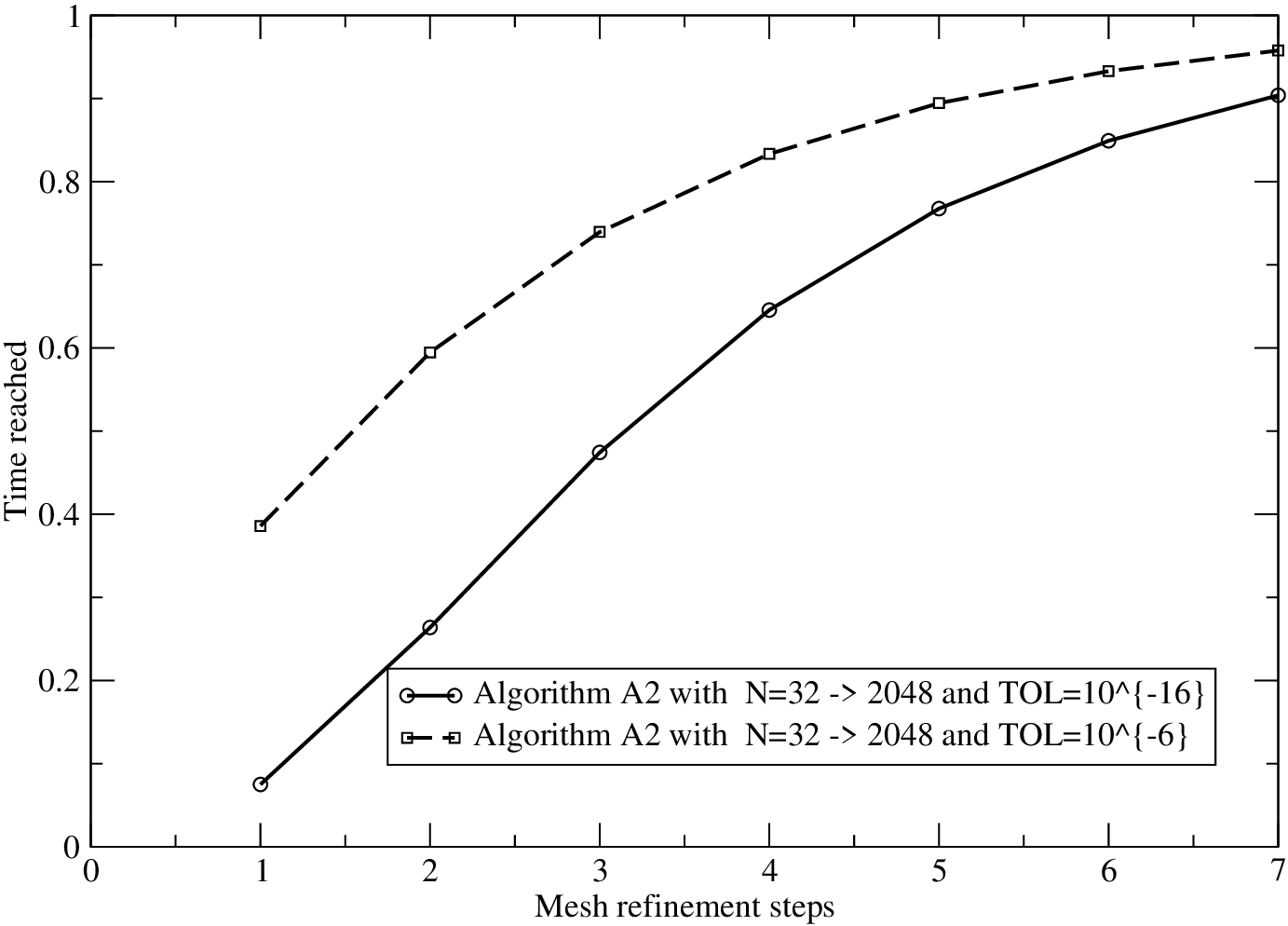,width=3.in}}
\caption{Case I, Galerkin, $\nu=0$: (a) Time spent between refinement steps for different tolerance values. (b) Time reached with the maximum allowed resolution. }
\label{plot_a2v0_galerkin}
\end{figure}

The results of the Algorithm 2 (A2) are presented in Figure \ref{plot_a2v0_galerkin}. We use this algorithm to perform mesh refinement when needed starting from a resolution $N_{start}=32$ and allowing a maximum resolution of $N_{final}=8192.$ We present results for two values of the tolerance $TOL1=10^{-16}$ and $TOL2=10^{-6}.$ We also plot the time reached with the maximum allowed resolution. When the tolerance criterion is less strict the algorithm can reach later times before running out of resolution. 

In Figure \ref{plot_a2v0_velocity} we compare the velocity field produced by A2 with $N_{start}=32,$ $N_{final}=2048$ and $TOL1=10^{-16}$ with the velocity field produced by A2 with $N_{start}=N_{final}=2048$ and the same tolerance. It is obvious that the results are in very good agreement. However, the mesh refinement calculation was about 13 times faster to run. The acceleration factor, with practically the same accuracy, can be increased by using a larger magnification ratio.

\begin{figure}
\centering
\epsfig{file=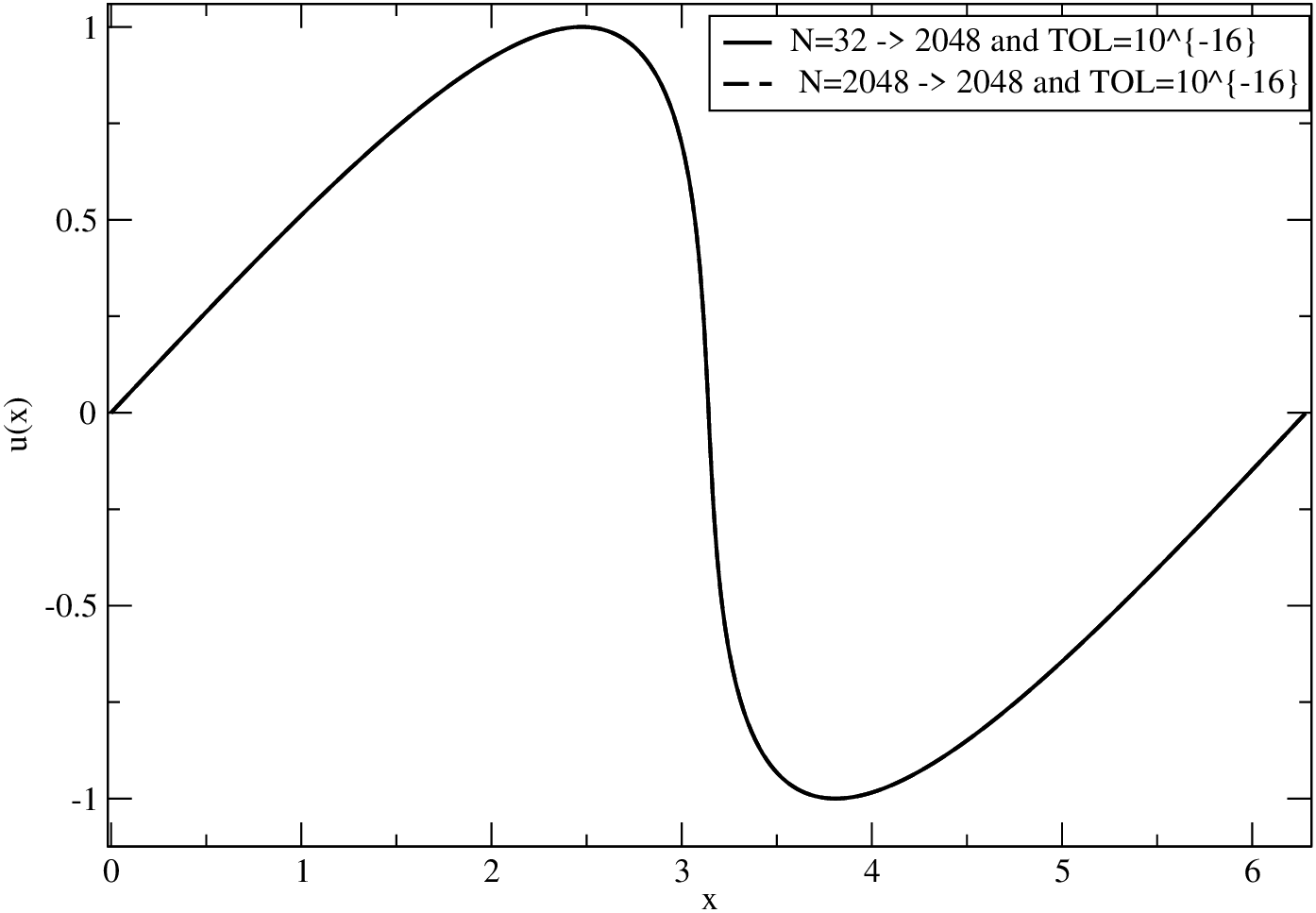,width=3.in}
\caption{Case I, Galerkin, $\nu=0$: Comparison of the velocity field produced at the time of termination of A2 for two different magnification ratios. 
The first simulation has $N_{start}=32$ and $N_{final}=2048$ while the second has $N_{start}=N_{final}=2048.$ }
\label{plot_a2v0_velocity_galerkin}
\end{figure}

Next, we present results for A3. The algorithm uses A2 until the maximum resolution is reached and then switches to the Galerkin model. We compare the energy evolution as predicted by A3 with $N_{start}=32, \; N_{final}=512$ and $TOL=10^{-16}$, a $t$-model calculation with $N=512$ and the random choice method \cite{chorinran} with $N=4096$ points. The mesh refinement algorithm produces results that are in very good agreement with the $t$-model and the random choice method {\it only} for times up to the formation of the singularity. Since the Galerkin model cannot dissipate energy it is bound to give erroneous predictions after the singularity is formed.    

\begin{figure}
\centering
{\epsfig{file=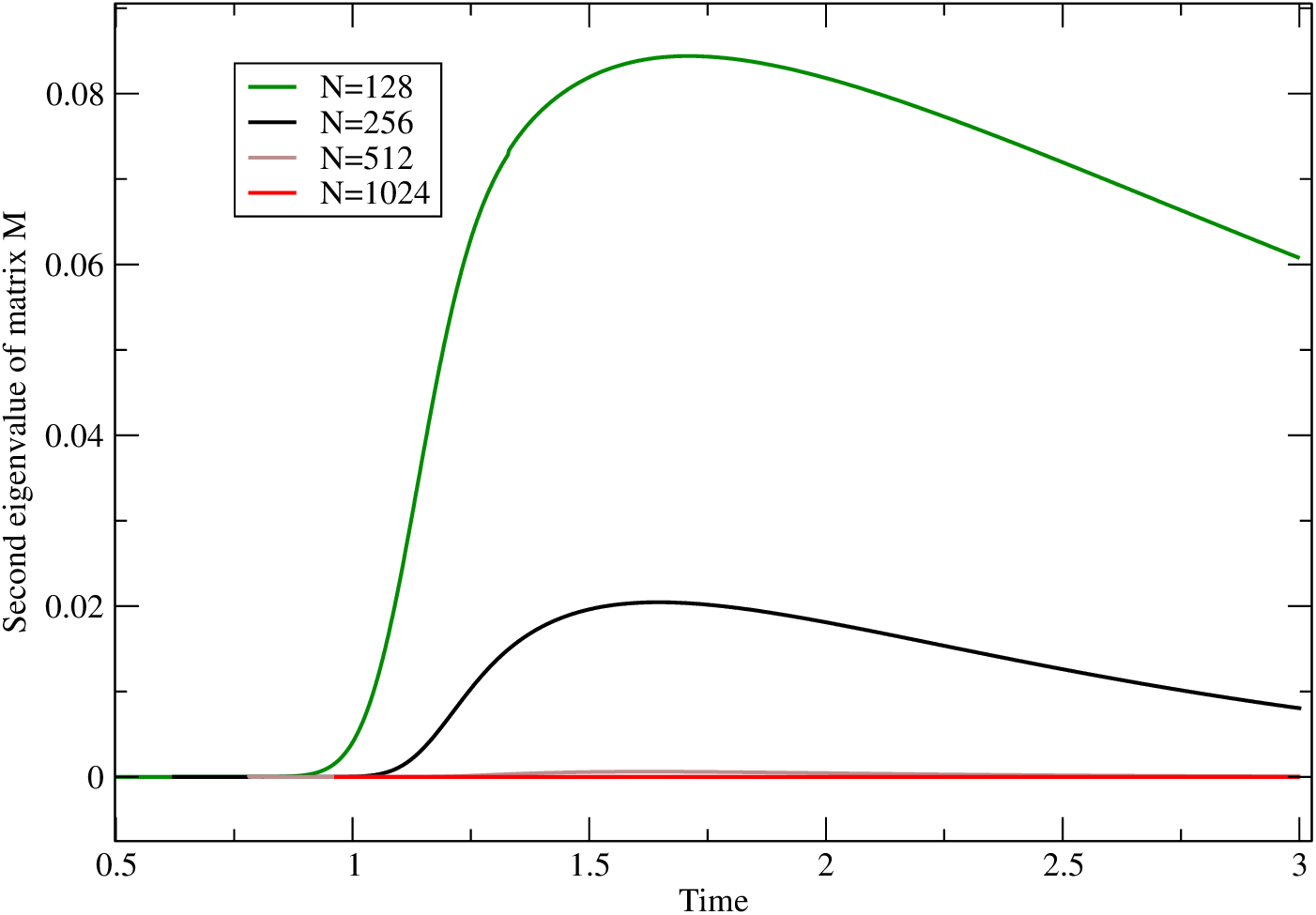,width=3.in}}
\caption{Case I, Galerkin, $\nu=0.01$: Evolution of the second eigenvalue of the matrix $M$ for different resolutions.}
\label{plot_a1v_2_galerkin}
\end{figure}

\begin{figure}
\centering
\subfigure[] {\epsfig{file=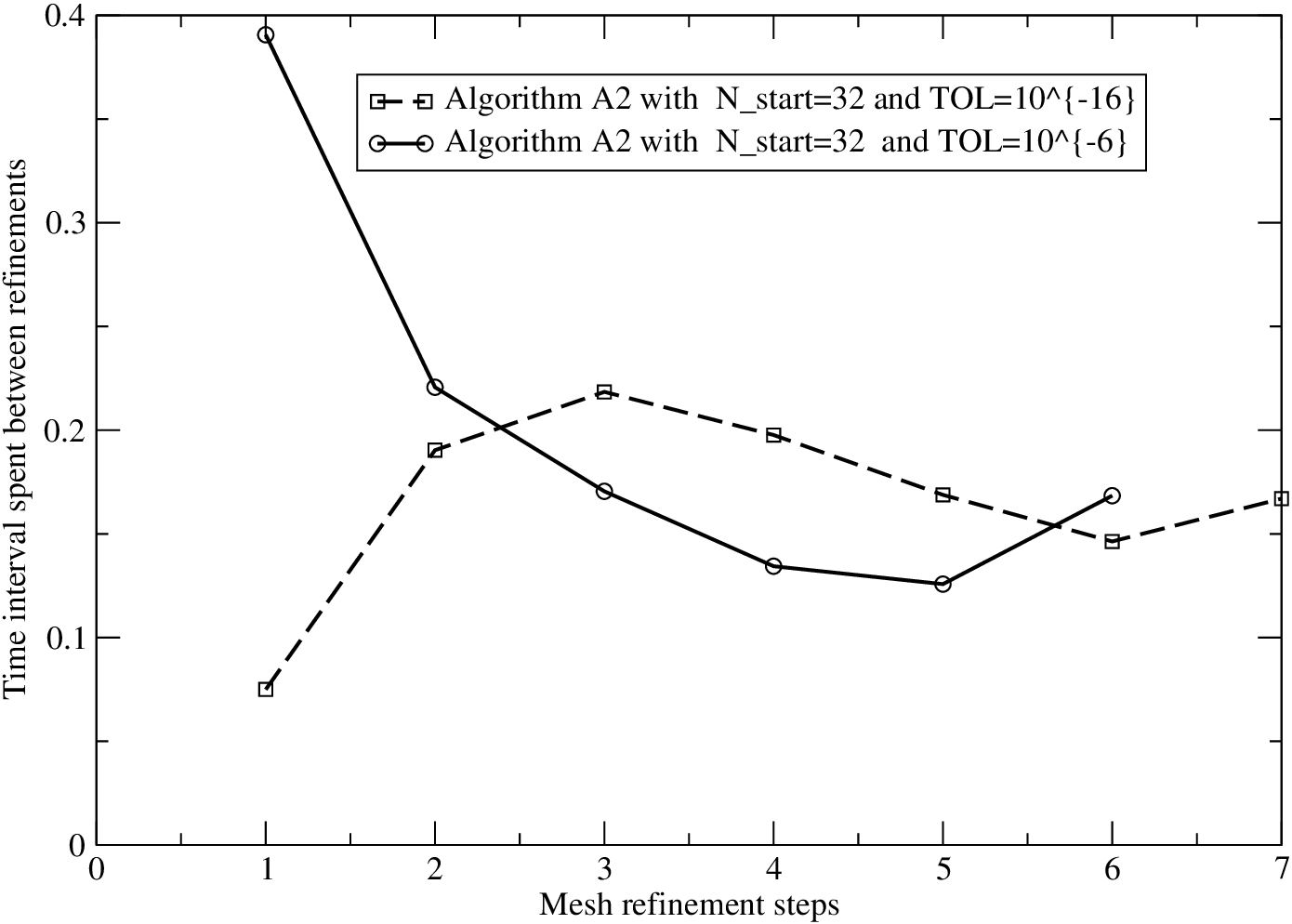,width=3.in}}
\qquad
\subfigure[] {\epsfig{file=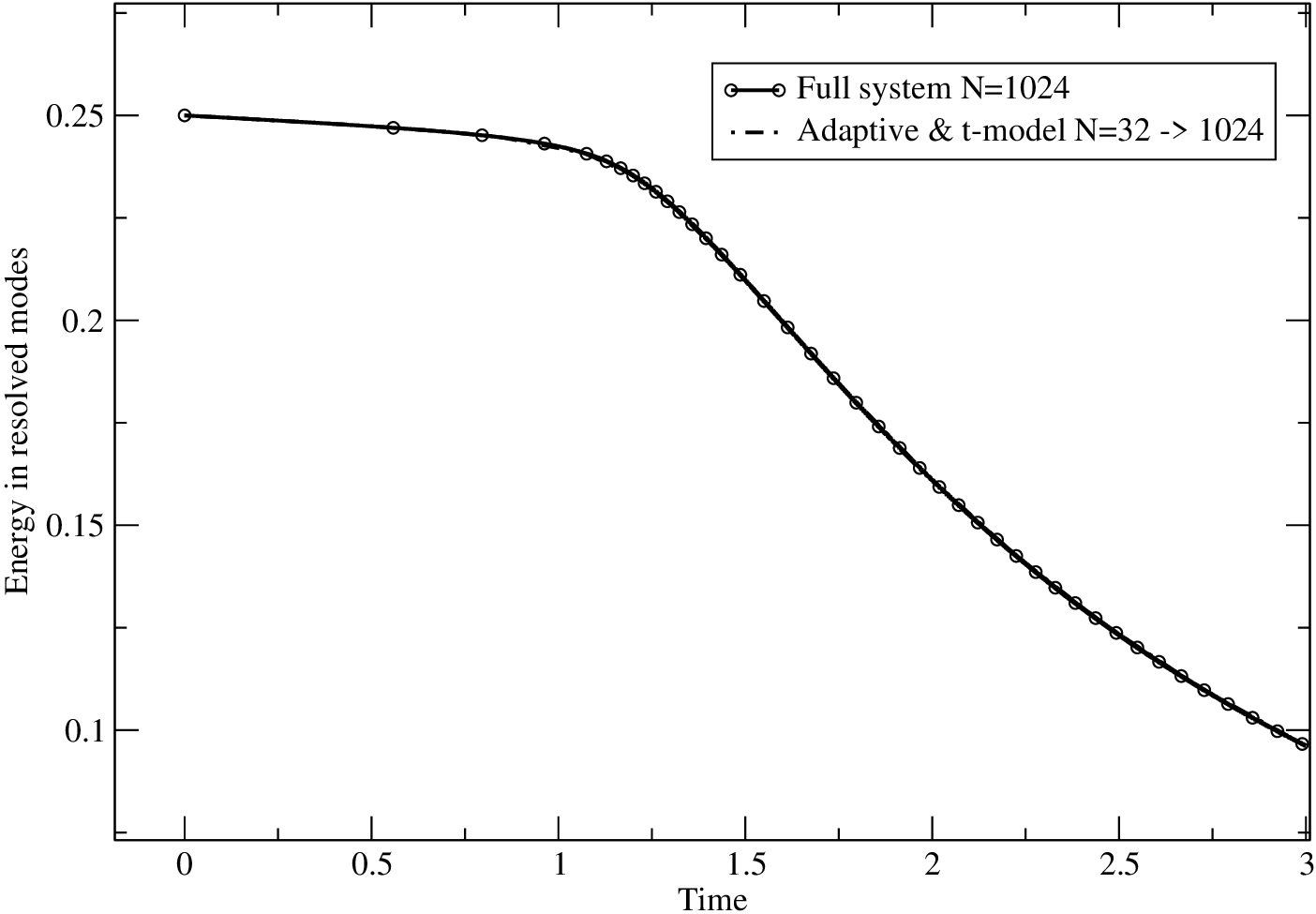,width=3.in}}
\caption{Case I, Galerkin, $\nu=0.01$: (a) Time spent between refinement steps for different tolerance values. (b) Comparison of the energy evolution as produced by A3 and the full system.}
\label{plot_a2v_2_galerkin}
\end{figure}

\paragraph{The case of nonzero viscosity} 

Figure \ref{plot_a1v_2_galerkin} shows the evolution of the second eigenvalue of the matrix $M$ for different resolutions with $\nu=0.01.$ Notice that the behavior of the eigenvalue is different from the case of no viscosity. There the eigenvalue increased with an increase in the resolution. Here it is decreasing until for large enough resolution it is practically zero for the duration of the simulation. In Figure \ref{plot_a2v_2_galerkin}(a) we present the progress of the refinement steps for $TOL1=10^{-16}$ and $TOL2=10^{-6}.$  We present the number of refinement steps needed until the refinement criterion is no longer activated. Figure \ref{plot_a2v_2_galerkin}(b) shows the evolution of the energy as predicted by A3 with $N_{start}=32, \; N_{final}=1024$ and $TOL=10^{-16}$ and a full system with $N=1024$ modes. The algorithm A3 produces results that are in very good agreement with the full system and does so 9 times faster.


\section{Case II: Knowledge of the terms but {\it not} of the coefficients of the reduced model}\label{know2}
We continue our presentation of the algorithms with the algorithms for the case when we have partial knowledge of the reduced model, i.e. knowledge of $\hR^{(1)}$ but {\it not} of $a^{(1)}.$ The constructions in this section can be considered as time-dependent generalizations of the Swendsen renormalization algorithm (e.g. see the nice presentation in Ch. 5 of \cite{binney}), even though here we do not have a statistical framework.

\subsection{How to locate a possible singularity}\label{location2}
When we do not have knowledge of the coefficients $a^{(1)},$ it is not possible to simulate the reduced model. In this case, the matrix $B$ has to be computed by using the values of the resolved modes from the full system. In this way the conditions \eqref{conditions} are automatically satisifed. However, note that the reduced and full systems can still deviate since the full system uses in addition the modes in $G$ while the reduced system does not. We can formulate an algorithm for locating the time of occurrence of a possible singularity.
\vskip14pt
{\bf Algorithm 4}
\begin{enumerate}
\item
Choose a resolution, i.e. a set of modes $F \cup G$ where $F$ are the resolved modes and $G$ the unresolved modes. Evolve the full system for the modes in $F \cup G.$ At the end of each step use the resolved modes as computed from the full system to compute the reduced system quantities. 
\item
Monitor the evolution of the eigenvalues of the renormalization matrix $M.$ 
\item
The eigenvalues corresponding to type i) terms should be 1 until the reduced system starts deviating from the full system. After an initial phase when the matrix $B$ is singular (to within the numerical precision used), the eigenvalues corresponding to type ii) terms should increase until they reach a maximum. Record the time evolution of the eigenvalues up to time of the maximum. This will guarantee that the time of occurrence of the turning point is included.
\item 
Repeat the experiment with a higher resolution. If the solution of the PDE develops a singularity at time $T^{*},$ the sequence (as we increase the resolution) of the instants of occurrence of the turning point should converge to $T^{*}.$ The converse is not necessarily true. Convergence of the time of occurrence of the eigenvalue turning point to a finite value, say $\tau,$ does not imply that the solution 
of the PDE develops a singularity at $\tau.$ It may well be that the solution remains smooth but requires higher resolution than currently available. 
\end{enumerate} 
The reader may be concerned that we are asking too much from our algorithm. We want to evolve only the full system and at the same time detect its progress towards underresolution through an internal check, i.e. without reference to an externally provided well resolved solution. This is not impossible. Since the reduced model calculations use only the resolved modes predicted by the full system, it is still possible to detect the deviation between the reduced and full models. While the full system is well resolved and energy is moving from $F$ to $G,$ the reduced and full systems will deviate at an increasing pace. When underresolution settles in, this deviation will slow down but it will still be present until aliasing effects dominate the energy specturm. Thus, a turning point will still be present close to the time of occurrence of a possible singularity. We are not advocating that this is a reliable method to monitor the deviation of the full and reduced system long after the occurrence of the singularity. In fact, no method that involves the full system can be accurate for long after the singularity occurs. As in Case I above, the more modes we use the more accurate the prediction of the turning point should be. In the limit  of an infinite number of modes, the turning point converges to the correct time of occurrence of the singularity.


\subsection{How to approach efficiently a possible singularity}\label{approach2}
The task of approaching efficiently a possible singularity can be accomplished through the same algorithm as in the Case I. The difference is that the matrix $B$ here is computed by using the resolved modes' values as computed from the full system.
\vskip14pt
{\bf Algorithm 5}
\begin{enumerate}
\item
Choose a value for $TOL.$ For this value of $TOL$ run a mesh refinement calculation, starting, say, from $N_{start}$ modes to $N_{final}$ modes. For example, let $N_{start}=32$ and double at each refinement until, say $N_{final}=256$ modes. Record the values of the quantities $\hat{E}_i, \; i=1,\ldots,m$ when $N=N_{final}$ and $|detB|=TOL.$ Let's call this simulation $S1.$
\item
For the same value of $TOL$ run a calculation with $N_{start}=N_{final}$ modes (for the example $N_{start}=N_{final}=256$). Record the values of the quantities $\hat{E}_i, \; i=1,\ldots,m$ when $|detB|=TOL.$ Let's call this simulation $S2.$
\item
Compare to within how many digits of accuracy the quantities $\hat{E}_i, \; i=1,\ldots,m$ computed from $S1$ and $S2$ agree. If the agreement is to within a specified accuracy, say 5 digits, then choose 
this value of $TOL.$ If the agreement is in fewer digits, then decrease $TOL$ (more stringent criterion) and repeat until agreement is met.
\item
Use the above decided value of $TOL$ to perform a mesh refinement calculation with a larger magnification ratio, i.e. a larger value for the ratio $N_{final}/ N_{start}.$  
\end{enumerate}
As in Case I, the agreement in digits of accuracy between $S1$ and $S2$ depends on the form of the terms chosen for the reduced model. Even though we do not know the coefficients of the reduced model, knowledge of the correct functional form of the terms can affect significantly the accuracy of the results. This situation is well known in the numerical study of critical exponents in equilibrium phase transitions (see e.g. Ch. 5 in \cite{binney}).


\subsection{How to follow a possible singularity}\label{follow2}
When we only know the functional form of the terms appearing in the reduced model but not their coefficients it is not possible to evolve a reduced system. Here we offer a partial fix of the problem which uses the same tools that we have been using so far to study the occurrence of a possible singularity. If the quantities $\hat{E}_i, \; i=1,\ldots,m$ are e.g. $L_p$ norms of the Fourier modes, then we can multiply Equations \eqref{reduced} with appropriate quantities and combine with Equations \eqref{conditions} to get
\begin{align*}\label{reduced}
\frac{d\hat{E}_1(\hat{u})}{dt} &= \sum_{i=1}^{m} a^{(1)}_i \hat{U}^{(1)}_{i1} (t,\hat{u}(t)) \\
\frac{d\hat{E}_2(\hat{u})}{dt} &= \sum_{i=1}^{m} a^{(1)}_i \hat{U}^{(1)}_{i2} (t,\hat{u}(t)) \\
 \quad     \cdots       \quad       & =   \quad       \cdots \quad \\
\frac{d\hat{E}_m(\hat{u})}{dt} &= \sum_{i=1}^{m} a^{(1)}_i \hat{U}^{(1)}_{im} (t,\hat{u}(t)) 
\end{align*}
where $\hat{U}^{(1)}_{ij}, \; i,j=1,\ldots,m$ are the new RHS functions that appear. Note that the RHS of the equations above does not involve primed quantities as in Case I. The reason is that here the reduced quantities are computed by using the values of the resolved modes from the full system. The above system of equations is a linear system of equations for the vector of coefficients $a^{(1)}.$ In fact, the matrix of the system is the transpose $B^T$ of the matrix $B.$ The linear system can be written as
\begin{equation}\label{alphasystem}
B^{T} a^{(1)}= {\bf e}
\end{equation}
where ${\bf e}=\bigl(\frac{d\hat{E}_1(\hat{u})}{dt}, \ldots, \frac{d\hat{E}_m(\hat{u})}{dt} \bigr).$ This system of equations can provide us with the time evolution of the vector $a^{(1)}.$ There is no need to compute the vector $a^{(1)}$ until the moment we reach the largest resolution available.  But, once the maximum resolution is reached, we can compute the vector $a^{(1)}$ and use it to define a reduced model which, hopefully, can follow the possible singularity. In summary, we have the  following algorithm:
\vskip14pt
{\bf Algorithm 6}
\begin{enumerate}
\item
Use {\bf Algorithm 5} until the maximum allowed resolution has been reached. Suppose that the set of modes for the maximum resolution 
is $K.$
\item
Divide $K$ in two sets $F$ and $G,$ where $F$ the resolved variables and $G$ the unresolved variables. Solve the system \eqref{alphasystem} to obtain the coefficient vector for the reduced model. 
\item
Use the reduced model computed at Step 2 to follow the modes in $F.$
\end{enumerate}
The determination of coefficients for the reduced model through the system \eqref{alphasystem} is a time-dependent version of the method of moments. We specify the coefficients of the reduced model so that the reduced model reproduces the rates of change of a finite number of moments of the solution. This construction can actually be used as an adaptive way of determining a reduced model if one has access to experimental values of the rates of change of a finite number of moments. Suppose that we are conducting a real world experiment where we can compute the values of a finite number of moments on a coarse grid only. Then we can use the system \eqref{alphasystem} at predetermined instants to update a model defined on the coarse grid. Results of this construction will be presented elsewhere.


\subsection{Numerical results for Case II}\label{numerics2}
As in Case I, we present results for the Burgers equation with zero and nonzero viscosity. Before we proceed with the details of the numerical results we want to make a general comment about the difference in actual computational time between the Case I and Case II algorithms. Recall that for Case I, we have to solve {\it both} the full and reduced system. The common stepsize is governed by the reduced model, since close to the singularity it starts losing energy and its stepsize decreases. This is just a result of the time derivative of the solution changing more rapidly close to the singularity. On the other hand, for Case II, we only solve the full system and thus there is no significant slowdown close to the singularity. Thus, for the same magnification ratio, the Case II algorithms turn out to be much faster than the Case I algorithms. However, the Case I algorithms have better convergence properties because they measure the deviation between the full system and an actual reduced model. This deviation is more pronounced for the same resolution and thus, for the same resolution the Case I algorithms give more accurate results. It would be interesting to conduct a detailed comparison of the two cases.

\subsubsection{The case of zero viscosity}\label{zero2}

\begin{figure}
\centering
\subfigure[]{\epsfig{file=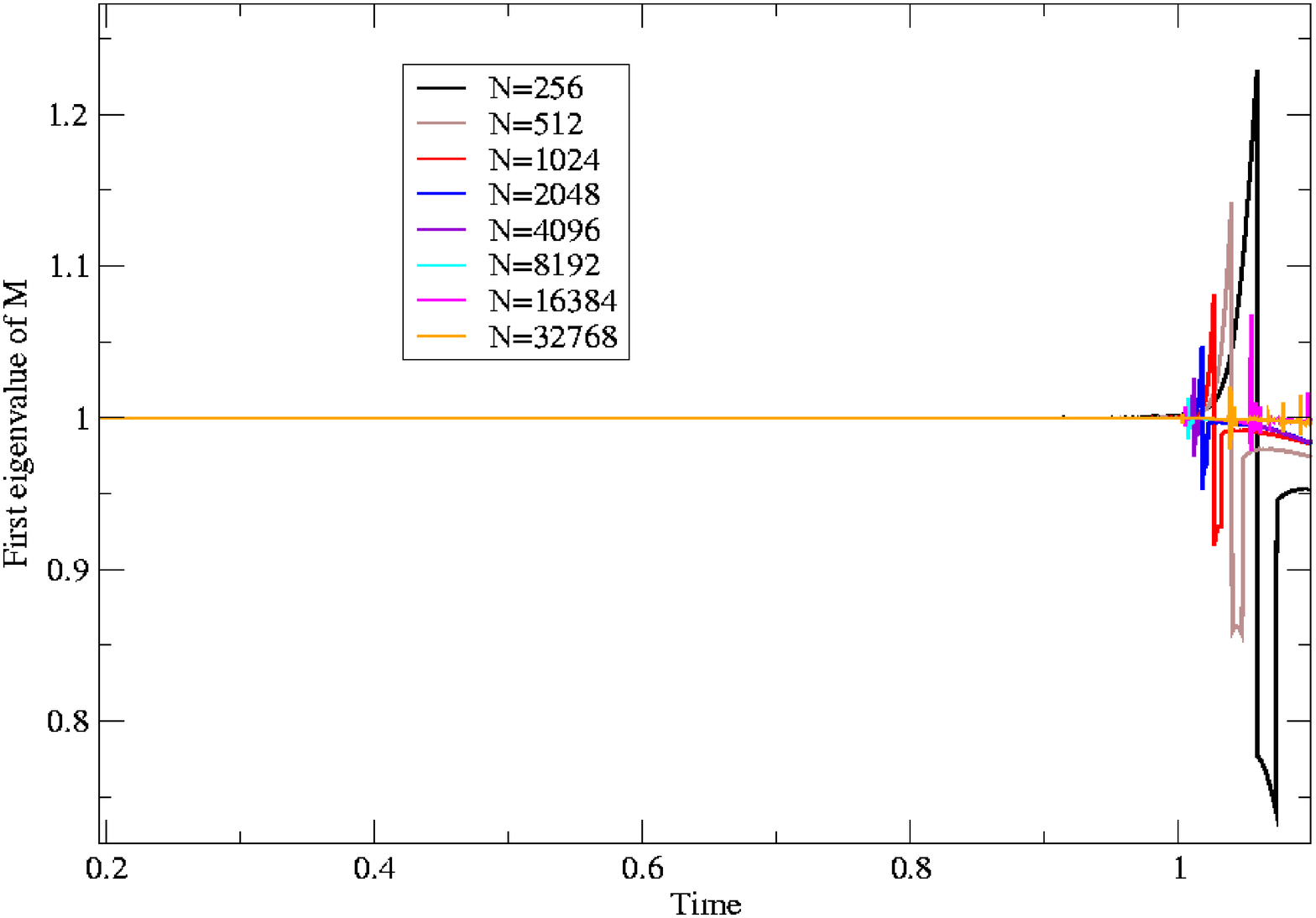, width=3.in}}
\qquad
\vskip14pt
\subfigure[]{\epsfig{file=plot_a4v0_b.eps,width=3.in}}
\caption{Case II, $\nu=0$: (a) Evolution of the first eigenvalue of the matrix $M$ for different resolutions. (b) Evolution of the second eigenvalue.}
\label{plot_a4v0}
\end{figure}

\begin{figure}
\centering
\epsfig{file=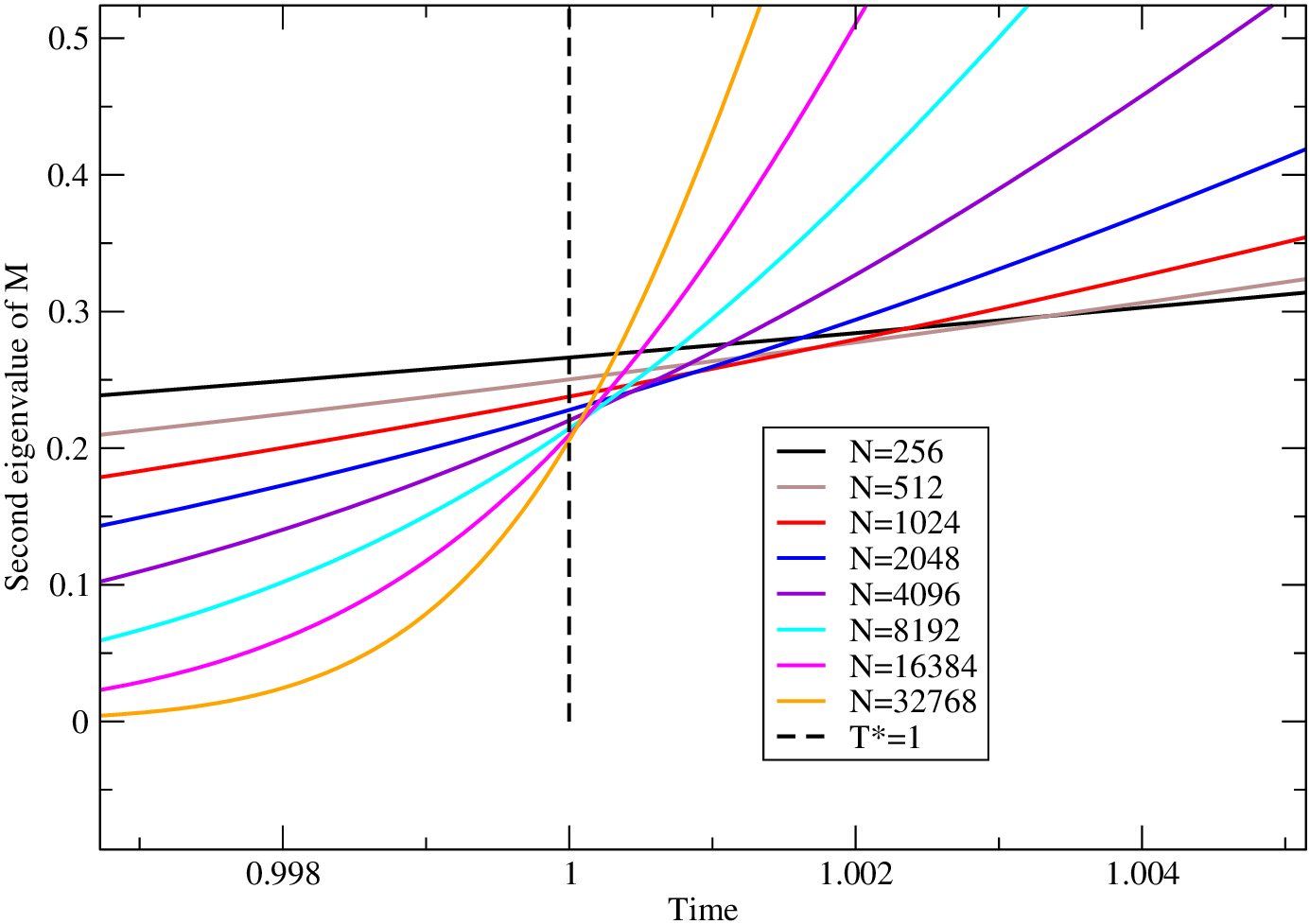,width=3.in}
\caption{Case II, $\nu=0$: Detail of the evolution of the second eigenvalue around the time of singularity. As the resolution is increased, the time when the evolution with a higher
resolution crosses the evolution with the lower resolution, moves towards 1.}
\label{plot_a4v0_detail_a}
\end{figure}

\begin{figure}
\centering
\epsfig{file=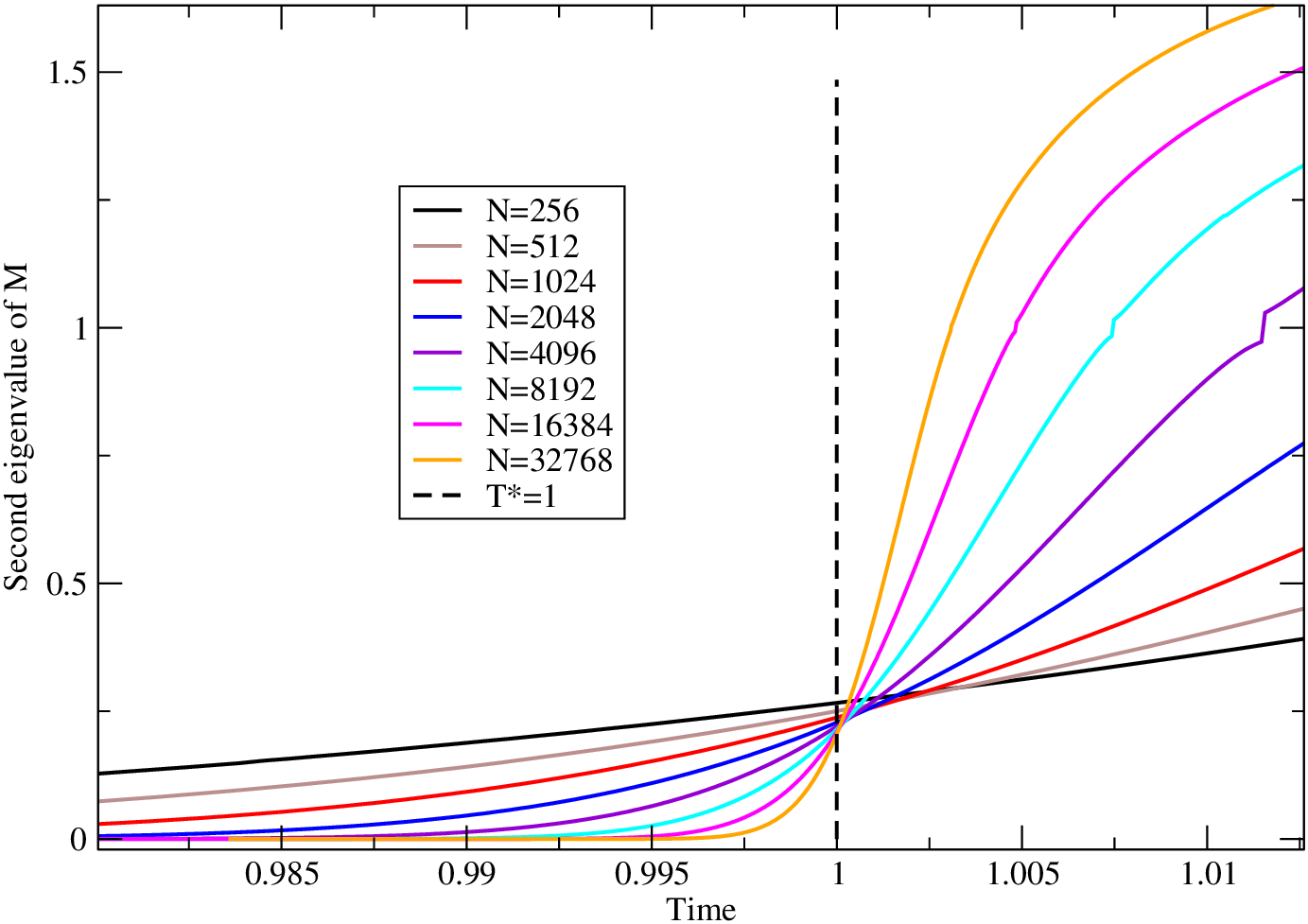,width=3.in}
\caption{Case II, $\nu=0$: Detail of the evolution of the second eigenvalue around the time of singularity. The turning point moves towards 1 as the resolution is increased.}
\label{plot_a4v0_detail_b}
\end{figure}

\begin{figure}
\centering
\subfigure[]{\epsfig{file=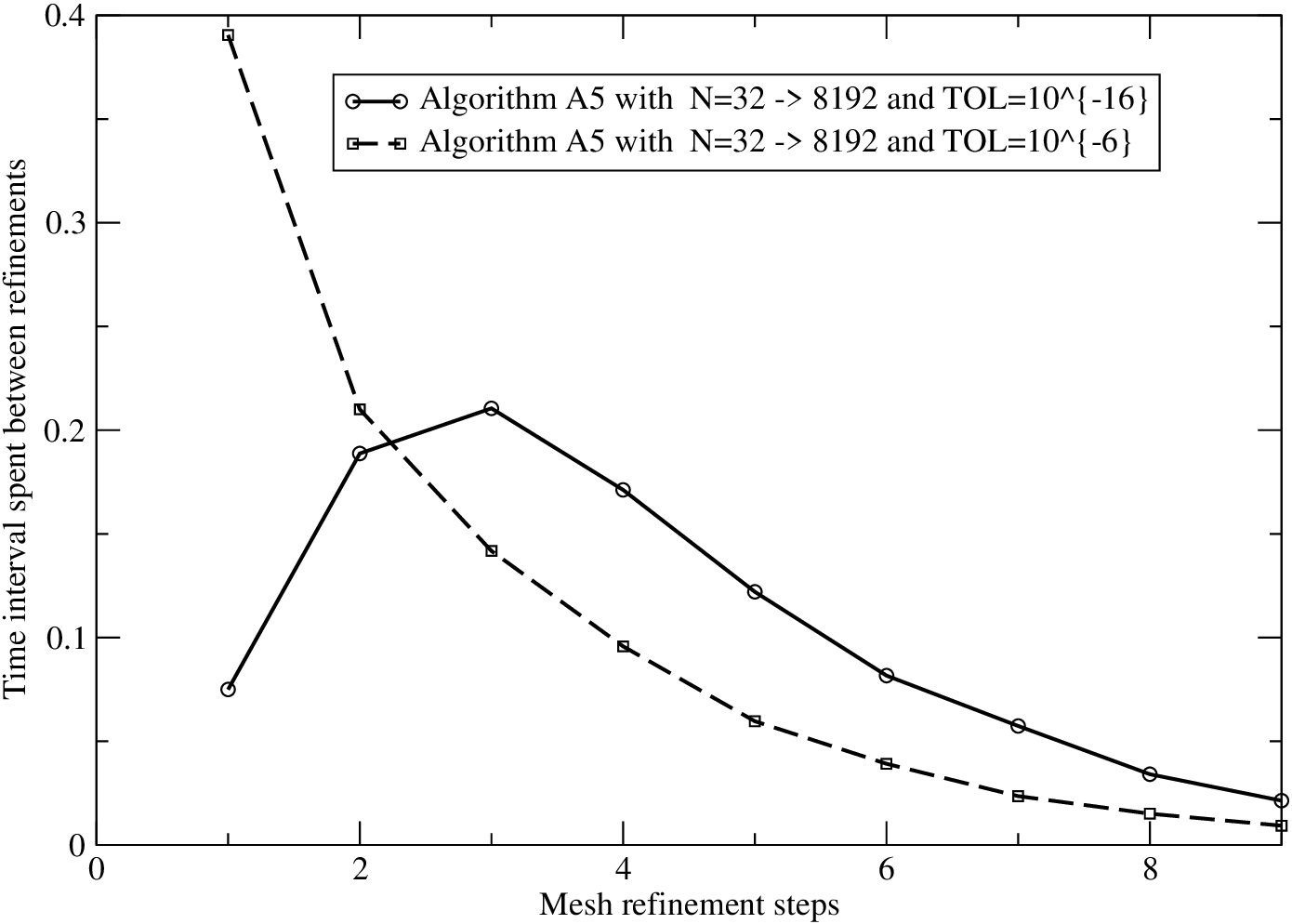,width=3.in}}
\qquad
\vskip20pt
\subfigure[]{\epsfig{file=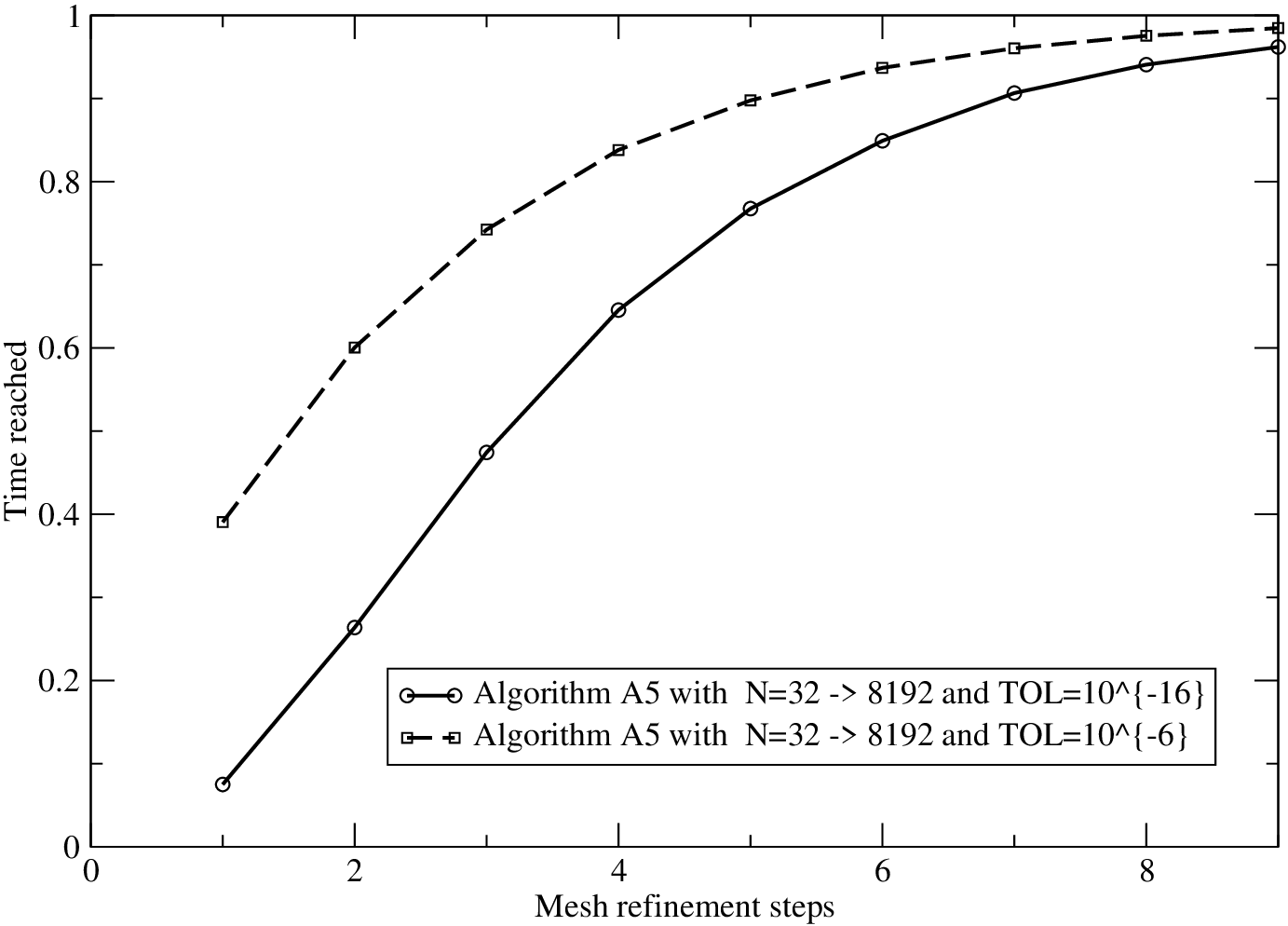,width=3.in}}
\caption{Case II, $\nu=0$: (a) Time spent between refinement steps for different tolerance values. (b) Time reached with the maximum allowed resolution. }
\label{plot_a5v0}
\end{figure}
                                                                                                                                                                                             
\begin{figure}
\centering
\epsfig{file=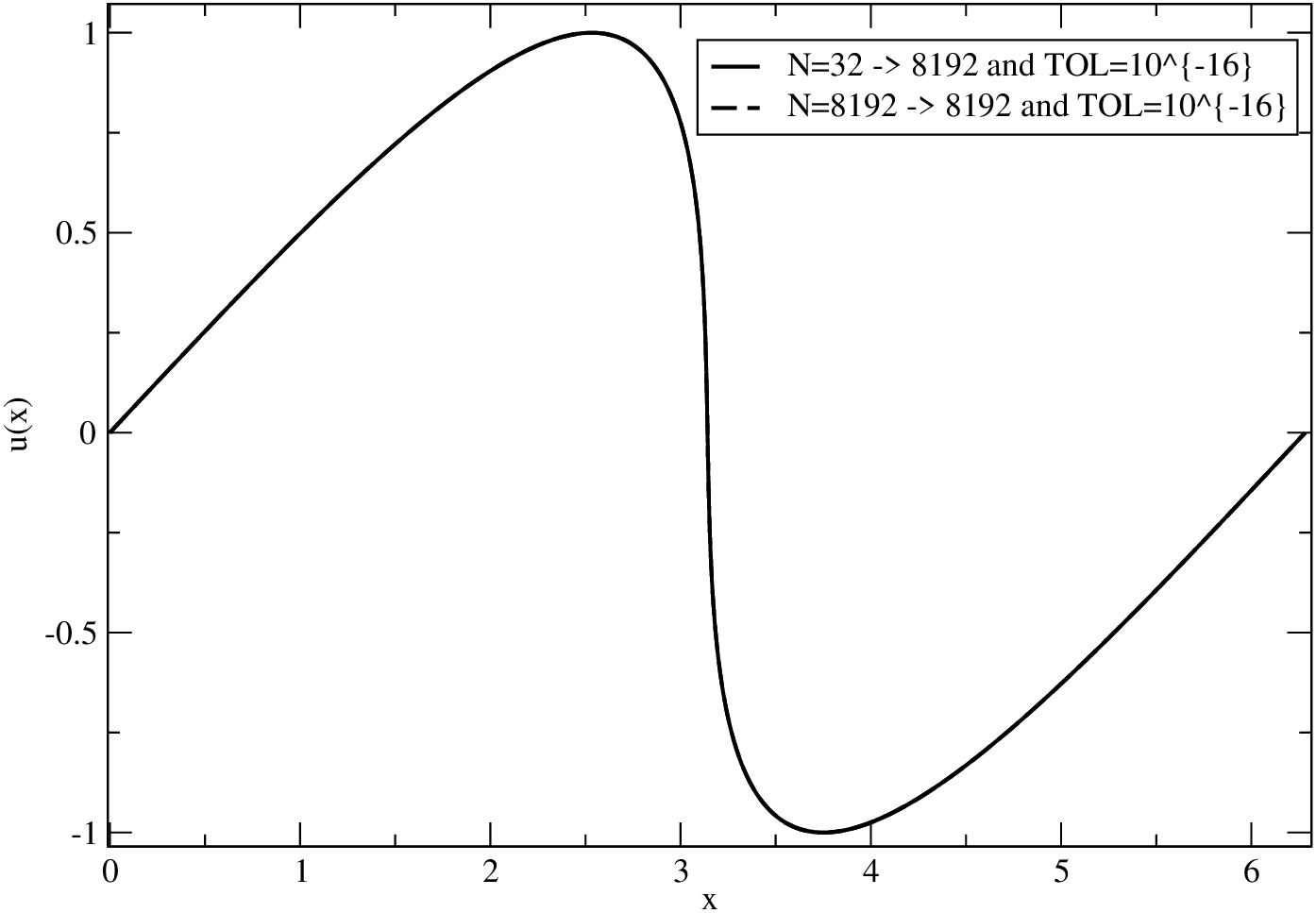,width=3.in}
\caption{Case II, $\nu=0$: Comparison of the velocity field produced at the time of termination of A5 for two different magnification ratios.
The first simulation has $N_{start}=32$ and $N_{final}=8192$ while the second has $N_{start}=N_{final}=8192.$ }
\label{plot_a5v0_velocity}
\end{figure}

\begin{figure}
\centering
{\epsfig{file=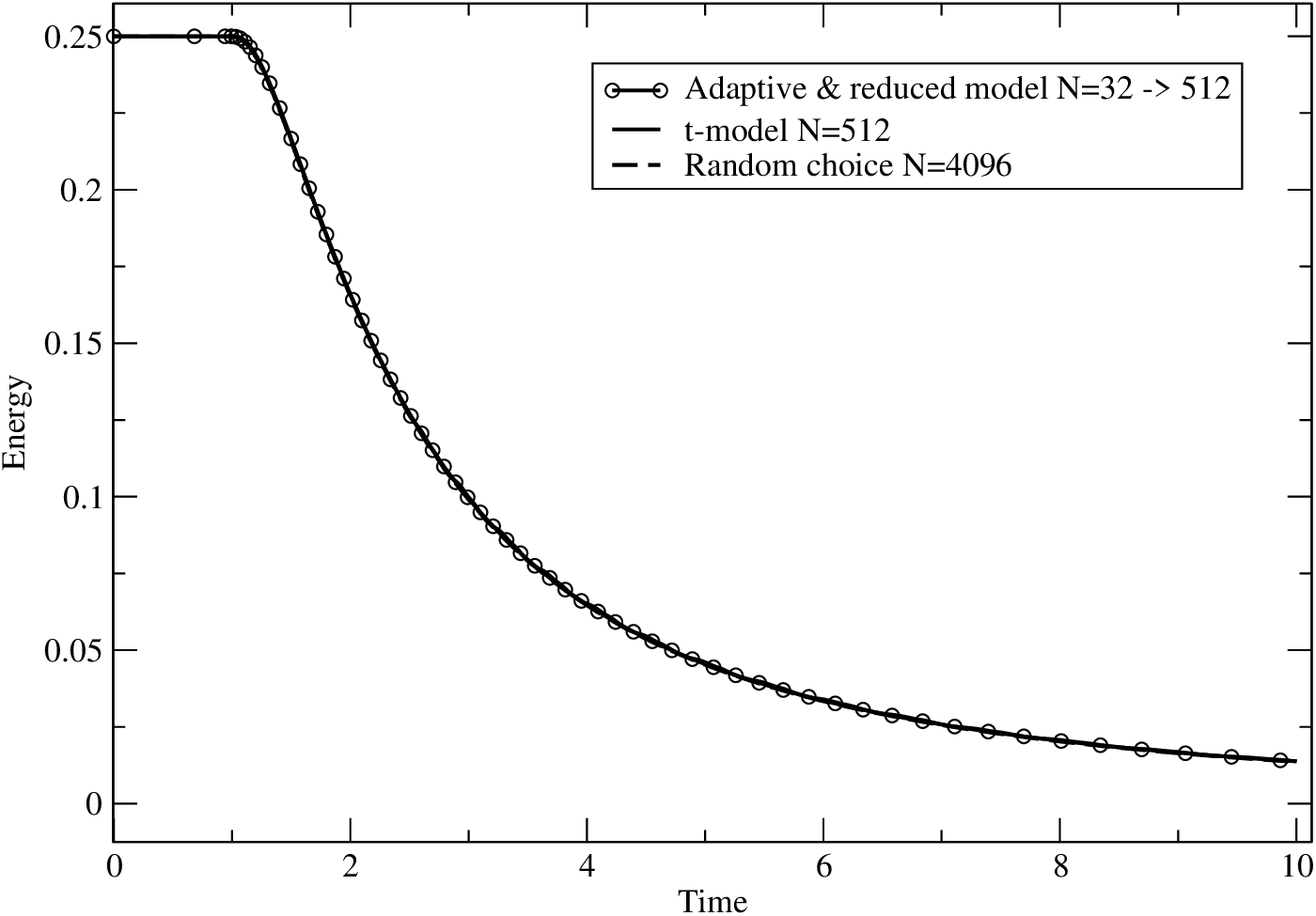,width=3.in}}
\caption{Case II, $\nu=0$: Comparison of the energy evolution as produced by A6, the $t$-model and the random choice method.}
\label{plot_a6v0_energy}
\end{figure}
                                                                                                                                                                                             
\begin{figure}
\centering
{\epsfig{file=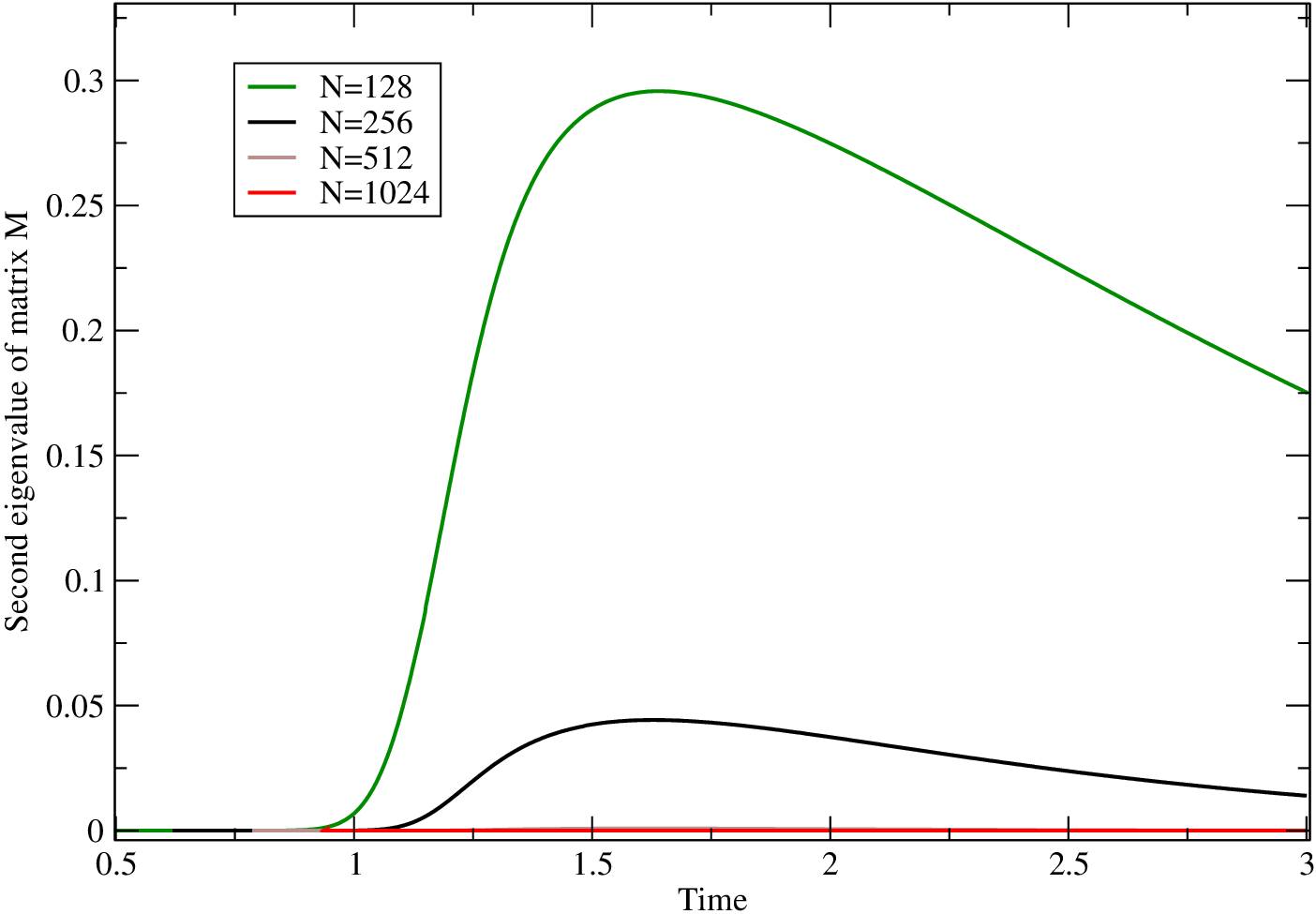,width=3.in}}
\caption{Case II, $\nu=0.01$: Evolution of the second eigenvalue of the matrix $M$ for different resolutions.}
\label{plot_a4v_2}
\end{figure}

\begin{figure}
\centering
\subfigure[] {\epsfig{file=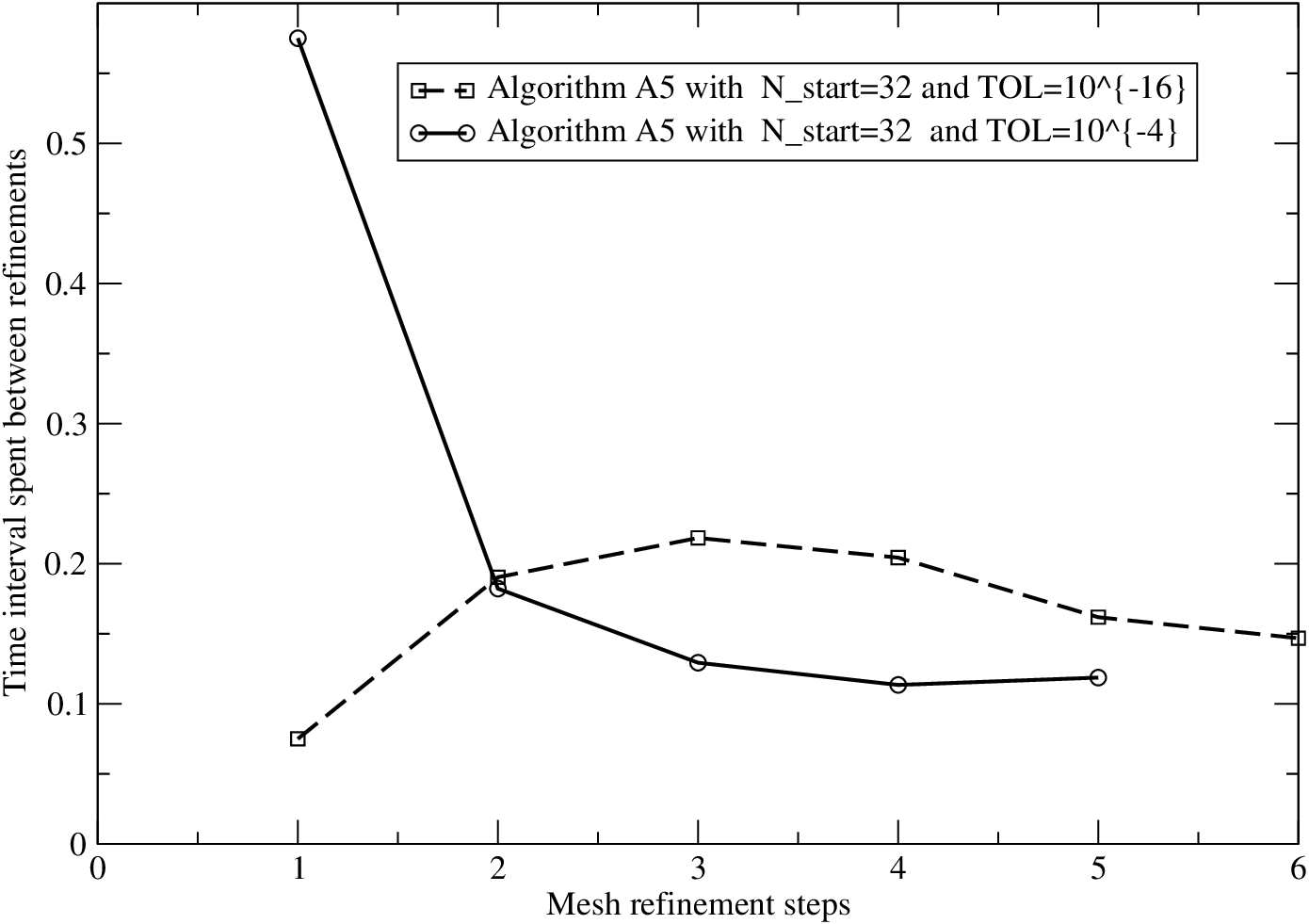,width=3.in}}
\qquad
\subfigure[] {\epsfig{file=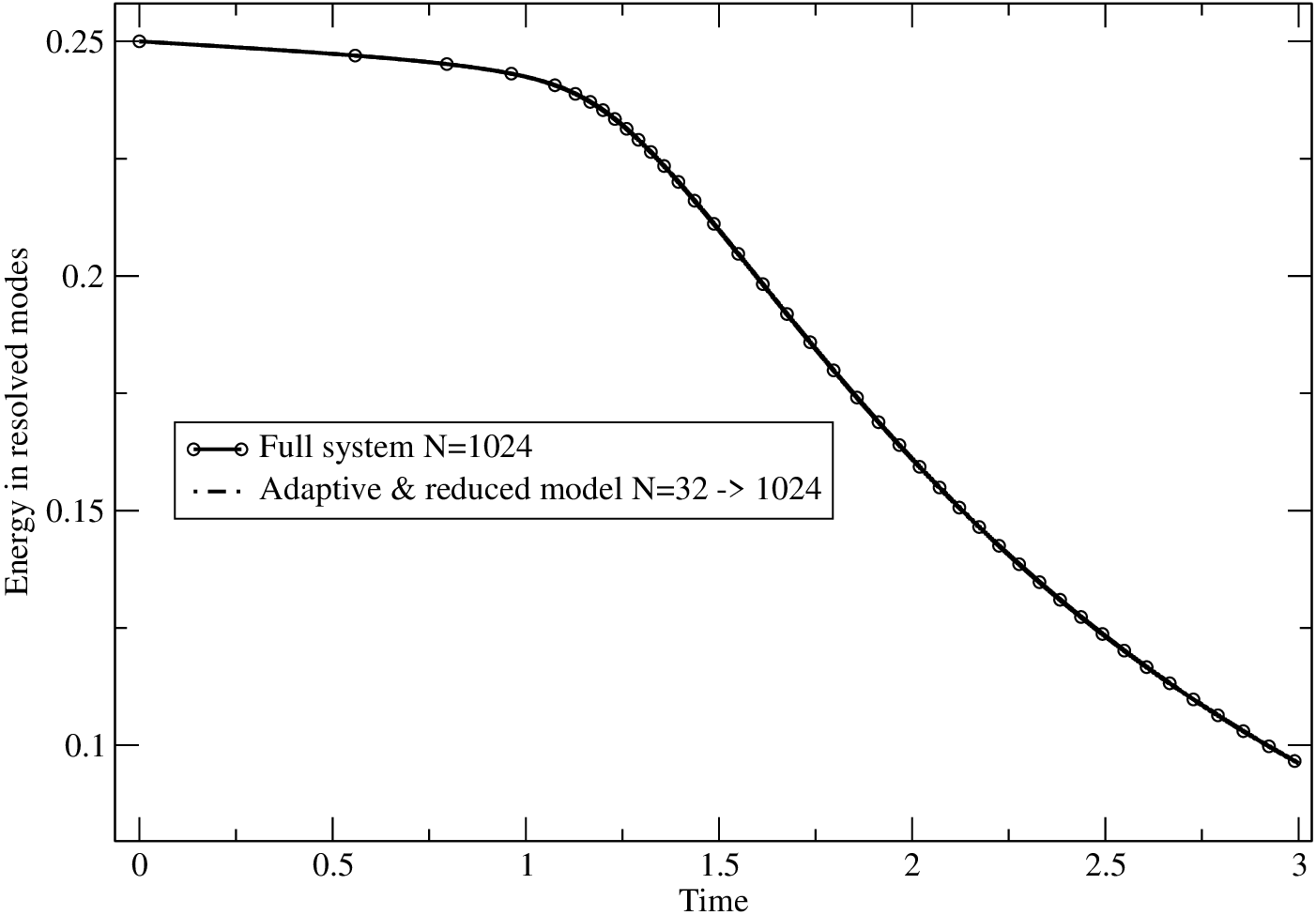,width=3.in}}
\caption{Case II, $\nu=0.01$: (a) Time spent between refinement steps for different tolerance values. (b) Comparison of the energy evolution as produced by A6 and the full system.}
\label{plot_a5v_2}
\end{figure}

First, we present results of the application of Algorithm 4 (from now on referred to as A4). Figure \ref{plot_a4v0} shows the evolution of the first and second eigenvalues of the matrix $M$ for different resolutions. As we expected, the first eigenvalue is equal to 1 until the moment when the reduced and full system start deviating. The second eigenvalue remains close to zero until the time of deviation of the reduced and full system. As in Case I, the evolution of the second eigenvalue is more informative about the behavior of the two systems. As we increase the resolution, the second eigenvalue remains closer to zero for a longer time. When it increases it does so at a rate that is higher the higher is the resolution. In Figures \ref{plot_a4v0_detail_a} and \ref{plot_a4v0_detail_b}  we have zoomed in on the time interval [.985,1.005]. As it should happen, the time when the evolution of the second eigenvalue for a higher resolution crosses the evolution of the second eigenvalue for the preceding resolution, moves towards 1. Also, the time of occurrence of the turning point moves towards 1 as we increase the resolution. The turning point for the highest resolution here $N=32768$ is around 1.002. For $N=8192,$ as the Case I calculation, the turning point is around 1.004.

The results of the Algorithm 5 (A5) are presented in Figure \ref{plot_a5v0}. We use this algorithm to perform mesh refinement when needed starting from a resolution $N_{start}=32$ and allowing a maximum resolution of $N_{final}=8192.$ We present results for two values of the tolerance $TOL1=10^{-16}$ and $TOL2=10^{-6}.$ As in Case I, when the tolerance criterion is less strict the algorithm can reach later times before running out of resolution. 

In Figure \ref{plot_a5v0_velocity} we compare the velocity field produced by A5 with $N_{start}=32,$ $N_{final}=8192$ and $TOL1=10^{-16}$ with the velocity field produced by A5 with $N_{start}=N_{final}=8192$ and the same tolerance. It is obvious that the results are in very good agreement. However, the mesh refinement calculation was about 240 times faster. For a calculation with $N_{start}=32, \; N_{final}=2048$ and the same tolerance (as the result presented for Case I) the acceleration factor was about 13.  

\begin{figure}
\centering
\subfigure[]{\epsfig{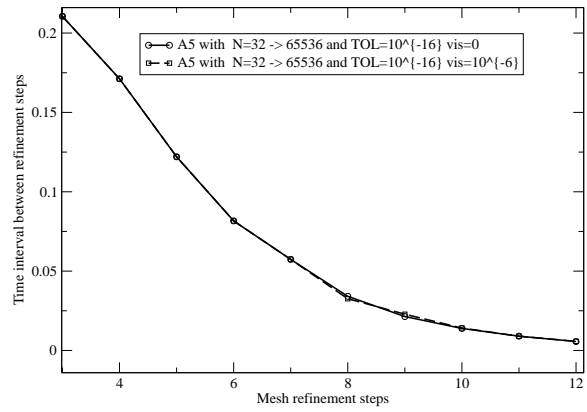}}
\qquad
\vskip20pt
\subfigure[]{\epsfig{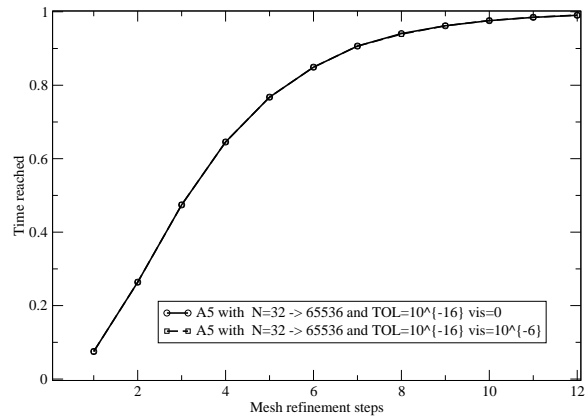}}
\caption{Case II, $\nu=0$ versus $\nu=10^{-6}$: (a) Time spent between refinement steps for two different viscosity values. (b) Time reached with the maximum allowed resolution for two different viscosity values. }
\label{plot_a5v0_compare}
\end{figure}

\begin{figure}
\centering
\epsfig{file=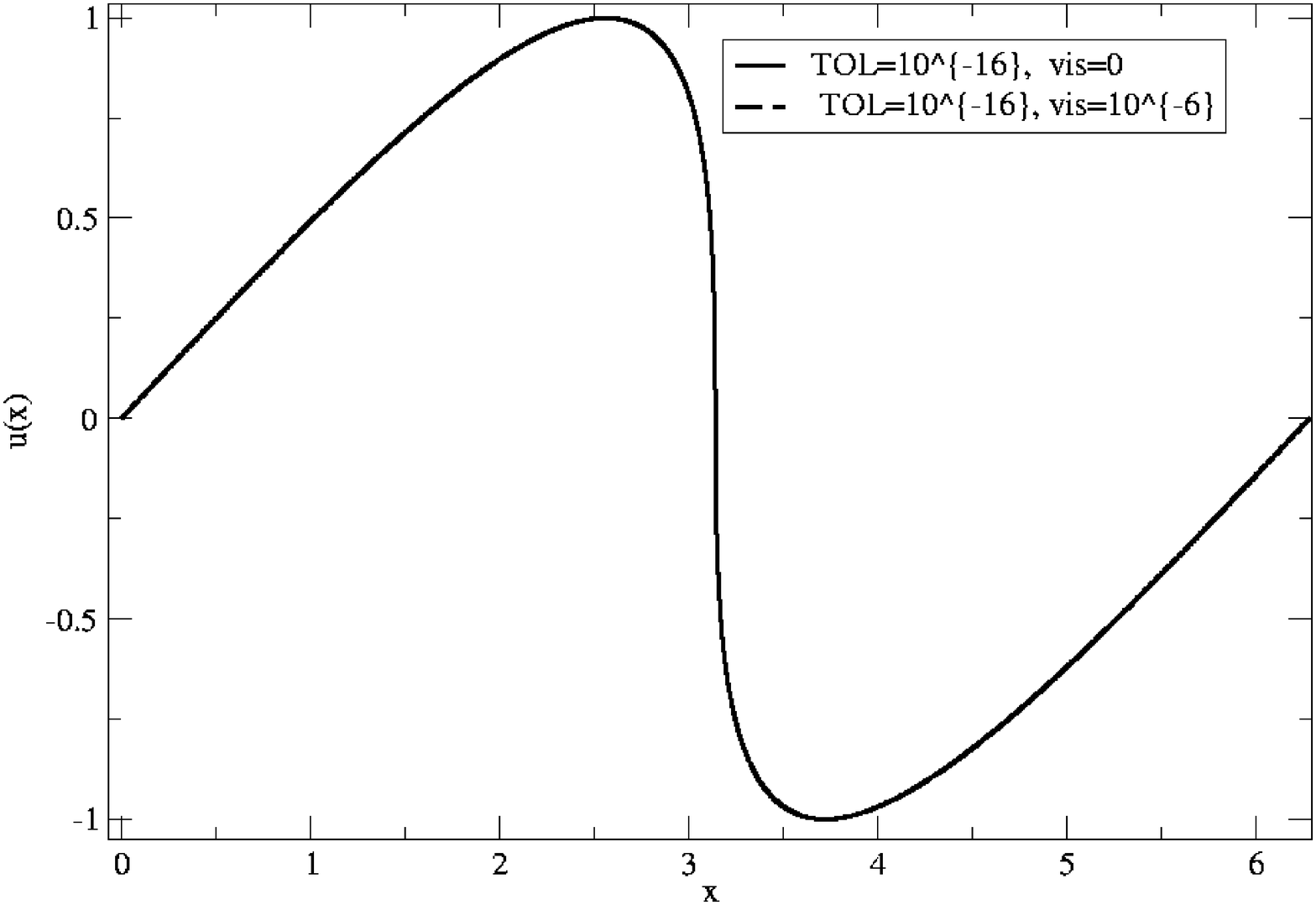,width=3.in}
\caption{Case II, $\nu=0$ versus $\nu=10^{-6}$: Comparison of the velocity field produced at the time of termination of A5 for two different viscosity values.}
\label{plot_a5v0_velocity_compare}
\end{figure}

Next, we present results for the Algorithm 6 (A6). The algorithm uses A5 until the maximum resolution is reached and then switches to the t-model. We compare the energy evolution as predicted by A6 with $N_{start}=32, \; N_{final}=512$ and $TOL=10^{-16}$, a $t$-model calculation with $N=512$ and the random choice method \cite{chorinran} with $N=4096$ points. When the maximum resolution is reached, the system \eqref{alphasystem} was solved once and the vector $a^{(1)}$ was determined. We find $a^{(1)}_1=1$ and $a^{(1)}_2=0.0284.$ The algorithm A6 produces results that are in very good agreement with the correct energy decay as computed by the random choice method.

\subsubsection{The case of nonzero viscosity}\label{nonzero2}
Figure \ref{plot_a4v_2} shows the evolution of the second eigenvalue of the matrix $M$ for different resolutions with $\nu=0.01.$ As in Case I, the second eigenvalue is decreasing until for large enough resolution it is practically zero for the duration of the simulation. In Figure \ref{plot_a5v_2}(a) we present the progress of the refinement steps for $TOL1=10^{-16}$ and $TOL2=10^{-4}.$  We present the number of refinement steps needed until the refinement criterion is no longer activated. In the first case, the criterion is very stringent and the algorithm has to refine more times to achieve it. Figure \ref{plot_a5v_2}(b) shows the evolution of the energy as predicted by A6 with $N_{start}=32, \; N_{final}=1024$ and $TOL=10^{-16}$ and a full system with $N=1024$ modes. The algorithm A6 produces results that are in very good agreement with the full system and does so 4 times faster.

\subsubsection{No viscosity versus very small viscosity}
We conclude our presentation with the comparison of the performance of the algorithms for the case of zero and finite but very small viscosity. Our aim is to show that a case with finite but very small viscosity cannot be distinguished from a case with zero viscosity if we do not have adequate resolution. Figures \ref{plot_a5v0_compare}(a) and (b) show the comparison of the results for the algorithm A5 with $N_{start}=32, \; N_{final}=65536$ and $TOL=10^{-16}$ for the case of $\nu=0$ and $\nu=10^{-6}.$ Figure \ref{plot_a5v0_velocity_compare} shows the velocity field at the moment when the algorithm runs out of resolution for the same magnification ratio and the two different viscosities. The two solutions are practically the same. Of course, if one allows higher resolutions eventually the viscous solution will deviate from the inviscid one.


\section{Conclusions and future work}\label{conclusions}
We have presented a collection of algorithms, inspired by renormalization constructions in critical phenomena, which allow the efficient location, approach and tracking of a possible singularity. The algorithms have two variants. The first one (Case I) assumes a complete knowledge of the reduced model. The second one (Case II) assumes knowledge of the functional form of the reduced model but not of the actual coefficients. The different algorithms were tested with rather encouraging results on the Burgers equation. On a theoretical level, the algorithms can be used to study the behavior of (near-) singular solutions. On the practical side, the algorithms can be used as a mesh refinement tool. We have only examined the simple case of periodic boundary conditions and the mesh refinement performed was uniform. An extension of the algorithms to a real space formulation is possible and we are currently working on it. A more detailed comparison of the Case I and Case II algorithms on other equations will be very interesting. Also, it would be interesting to investigate more the adaptive reduced model construction as described at the end of Section \ref{follow2}.

The problem of how to construct mesh refinement methods and how to approach more efficiently a possible singularity has attracted considerable attention (see e.g. \cite{almgren,berger1,berger2,budd,ceniceros,landman}). We plan to extend the constructions presented here to a real space formulation which will allow the treatment of non-periodic boundary conditions and more complicated geometries. Also, one can think of generalizing the algorithms to allow for mesh coarsening after the phase of the interesting, and computationally intensive, dynamics has passed.

The original motivation behind the development of the algorithms was the open problem of the formation of singularities in finite time for the Euler and Navier-Stokes equations of fluid mechanics. As we have stressed, the algorithms allow us to come close to a possible singularity but not decide definitely on the presence of an actual singularity. Maybe the algorithms can be used to decide such an issue although it is not clear to the author at the current stage.  However, we hope that the algorithms can be of use to the simulation of real world flows by allowing a better assessment of the onset of underresolution.

\section*{Acknowledgements} I am grateful to Profs G.I. Barenblatt, 
A.J. Chorin, O.H. Hald and Dr. Y. Shvets for many helpful discussions and comments. This work was supported in part by the Director, Office of Science, Advanced Scientific Computing Research, U.S. Department of Energy under Contract No. DE-AC02-05CH11231.


\begin{thebibliography}{99}


\bibitem{almgren}
Almgren A.S., Bell J.B., Colella P., Howell L.H. and Welcome M.L., A conservative adaptive projection method for the variable-density incompressible Navier-Stokes equations, J. Comp. Phys. 142 (1998) pp. 1-46.

\bibitem{berger1}
Berger, M. and Kohn, R., A rescaling algorithm for the numerical calculation of blowing-up solutions, 
Comm. Pure Appl. Math. 41 (1988) pp. 841-863.


\bibitem{berger2}
Berger M. and Colella P., Local adaptive mesh refinement for shock hydrodynamics, J. Comp. Phys. 82 (1989) pp. 62-84.

\bibitem{bernstein}
Bernstein D., Optimal prediction of Burger's equation, Multi. Mod. Sim. 6 (2007) pp. 27-52.


\bibitem{binney}
Binney J., Dowrick N., Fisher A., Newman M., The Theory of Critical Phenomena (An Introduction to 
the Renormalization Group), The Clarendon Press, Oxford, 1992.

\bibitem{boyd}
Boyd J.P., Chebyshev and Fourier Spectral Methods, Dover, New York, 2001.

\bibitem{budd}
Budd, C. J., Huang, W. and Russell, R. D., Moving mesh methods for problems with
blow-up. SIAM Jour. Sci. Comput. 17 (1996) pp. 305-327.

\bibitem{ceniceros}
Ceniceros H.D. and Hou T.Y., An efficient dynamically adaptive mesh for potentially singular solutions, J. Comp. Phys. 172 (2001) pp. 609-639.

\bibitem{chorinran}
Chorin A.J., Random choice solution of hyperbolic systems, J. Comp. Phys. 22 (1976) pp. 517-533.

\bibitem{CHK00}
Chorin, A.J., Hald, O.H. and Kupferman, R.,
Optimal prediction and the Mori-Zwanzig representation of irreversible
processes. Proc. Nat. Acad. Sci. USA 97 (2000) pp. 2968-2973.

\bibitem{CHK3}
Chorin, A.J., Hald, O.H. and Kupferman, R., Optimal prediction with memory, Physica D 166 (2002) pp. 239-257.

\bibitem{CS05}
Chorin, A.J. and Stinis, P., Problem reduction, renormalization and memory,
Comm. App. Math. Comp. Sci. 1 (2005) pp. 1-27.

\bibitem{goldenfeld}
Goldenfeld, N., Lectures on Phase Transitions and the Renormalization Group, Perseus Books, Reading, Mass., 1992.

\bibitem{givon}
Givon, D., Kupferman, R. and Stuart, A., Extracting macroscopic 
dynamics: model problems and algorithms, Nonlinearity 17 (2004) pp. R55-R127.

\bibitem{GSW}
Givon, D., Stinis, P. and Weare, J., Dimensional reduction for particle filters of systems with time-scale separation Technical Report LBNL-62141 (2007) (submitted to IEEE Trans. Signal Proc.).

\bibitem{hair}
Hairer, E.,  N\"orsett, S.E.  and Wanner, G., Solving Ordinary Differential Equations I, Springer, NY, 1987.

\bibitem{HS06}
Hald O.H. and Stinis P., Optimal prediction and the rate of decay for solutions of the Euler equations in two and three dimensions, Proc. Natl. Acad. Sci.,104, no. 16 (2007) pp. 6527-6532.

\bibitem{landman}
Landman, M.J.,  Papanicolaou, G.C., Sulem, C. and Sulem, P.,  Rate of blowup for solutions of the nonlinear Schršdinger equation at critical dimension, Phys. Rev. A 38 (1988) pp. 3837-3847.

\bibitem{S05}
Stinis P., A maximum likelihood algorithm for the estimation and renormalization of exponential densities, J. Comp. Phys. 208 (2005) pp. 691-703. 

\bibitem{weinberg}
Weinberg S., Why the renormalization group is a good thing, Asymptotic Realms of Physics: Essays in Honor of Francis E. Low, MIT Press, Cambridge MA, 1983.



\end{thebibliography}
\end{document}